\documentclass[a4paper, 12pt]{article}
\usepackage[dvips]{graphicx}
\usepackage{amsfonts}
\usepackage{amssymb}
\usepackage{color}
 
\def\a{{\alpha}}
\def\B{{\cal B}}
\def\b{{\beta}}
\def\C{{\cal C}}
\def\g{{\gamma}}

\def\d{{\delta}}
\def\E{{\cal E}}
\def\e{{\epsilon}}
\def\F{{\cal F}}

\def\l{{\lambda}}
\def\L{{\cal L}}
 \def\P{{\cal P}}
 \def\Q{{\cal Q}}
\def\s{{\sigma}}
\def\S{{\cal S}}
\def\t{{\theta}}

\def\CC{{\mathbb C}}
\def\Chat{\hat{\mathbb C}}

\def\HH{{\mathbb H}}
\def\NN{{\mathbb N}}
\def\RR{{\mathbb R}}
\def\QQ{{\mathbb Q}}
\def\ZZ{{\mathbb Z}}

\def\ch {{\cosh}}
\def\sh {{\sinh}}
\def\th {{\tanh}}
\def\cth {{\coth}}

\def\dd{ {\partial}}
\def\ax{\mathop{ Ax}}
\def\bch{{\partial{\cal C}}}
\def\Int{\mathop{\rm  Int}}
\def\ML{{\cal  ML}}
\def\PML{{\rm P}{\cal {ML}}}
\def\pl{ pl}
\def\QF{{\cal Q}{\cal F}}
\def\Sf{ S}
\def\teich{ {\cal F}}
\def\Teich{{\cal F}}
\def\tr{\mathop{\rm Tr}}
\def\p_#1{ \frac{\partial \hfil}{\partial {#1}}}

\newtheorem{thm}{Theorem}[section]
\newtheorem{lemma}[thm]{Lemma}

\newtheorem{cor}[thm]{Corollary}
\newtheorem{prop}[thm]{Proposition}

\newtheorem{remark}[thm]{Remark}
\newenvironment{proof}{{\sc Proof.}}{$\;\square$ \vskip .2in}

\def\qed{{$\;\square$ \vskip .2in}}

\begin{document}

\title{Limits of quasifuchsian groups
with small bending}
\author{Caroline
Series\\Mathematics Institute, Warwick University, \\Coventry CV4
7AL,
U.K.}

\maketitle

\bibliographystyle{Plain}

\begin{abstract} We study limits of quasifuchsian groups
for which the bending measures on the convex hull boundary  tend to
zero, giving  necessary and sufficient conditions for the limit group
to exist and be Fuchsian. 
As an application we complete the proof of a conjecture made
in~\cite{S1},  that the closure of pleating varieties for
quasifuchsian groups meet Fuchsian space exactly in  Kerckhoff's lines
of minima of length functions. Doubling our examples gives rise to a
large class of cone manifolds which degenerate to hyperbolic surfaces as
the cone angles
approach $2\pi$.

\medskip
\noindent {\small AMS classification numbers: 30F40, 20H10, 32G15.}
\end{abstract}

\section{Introduction}
\label{sec:introduction}

Considerable recent interest has focussed on the two components of the
convex hull boundary of a quasifuchsian group $G$. The object of this
paper is to study what happens when these components flatten out,  the
obvious expectation being that under
suitable conditions a limit group should exist and be Fuchsian.

 The   hyperbolic $3$-manifold $\HH^3/G$ associated to the group $G$
is homeomorphic to $\Sf \times (0,1)$ for some topological surface
$\Sf$. The convex hull boundary, that is, the boundary of the convex
hull of all closed geodesics in $\HH^3/G$,  
has two connected components each themselves homeomorphic to  $\Sf$.
Each
component is  bent along some geodesic
lamination on $\Sf$, the amount of bending being measured by the {\it
bending measures}  $\pl^{\pm} = \pl^{\pm}(G)$. 
Given two
measured laminations $\mu$ and $\nu$,  the \emph{pleating
variety} $\P_{\mu,\nu}$ consists of all groups $G \in \QF(\Sf)$ for
which $pl^+(G)$ is projectively equivalent to $\mu$ and 
 $pl^-(G) $ to $\nu$.

 Recall that measured laminations $\mu$ and $\nu$
are said to  fill  up $\Sf$ if $i(\mu,\xi) + i(\nu,\xi)>0$ for any
measured lamination $\xi$. It is not hard to show that  $\P_{\mu,\nu}$
is empty unless  $\mu$ and $\nu$ fill  up $\Sf$.  
The converse is a  special case of central recent result
of Bonahon-Otal~\cite{BonO}:
\begin{thm}
\label{thm:otalbon} Let $\Sf$ be a hyperbolic surface and let $\mu,\nu$
be measured geodesic laminations which fill 
up $\Sf$. Then there is a quasifuchsian group $G(\mu,\nu)$ for which 
$pl^+(G) =\mu$ and $pl^-(G)=\nu$. If $\mu,\nu$ are rational, then 
$G(\mu,\nu)$ is unique. 
\end{thm}
  One could well conjecture that the final uniqueness statement  is true
in general. Thus from now on we use $G(\mu,\nu)$ to denote any
quasifuchsian group
for which $pl^+(G) =\mu$ and $pl^-(G)=\nu$. We prove:
\begin{thm}
\label{thm:main} Let $\mu, \nu $ be 
two measured laminations which together fill up $\Sf$.
 Then as $\t \to 0$, the sequence
$G(\theta \mu, \theta \nu )$  converges to a Fuchsian group.
\end{thm}
 (The Bonahon-Otal result for irrational
laminations involves a delicate limit process which however says nothing
about what happens when the bending measures tend to zero.)

We also identify the limit Fuchsian group precisely.
 Based on Thurston's earthquake theorem, in~\cite{KerckLM}, Kerckhoff
proved the following result about length functions on Teichm\"uller
space:
\begin{thm}
\label{thm:kerck} Let $\Sf$ be a hyperbolic surface and let $\mu,\nu$
be measured geodesic laminations which fill 
up $\Sf$. Then the length function $l_{\mu} + l_{\nu}$ has a unique
minimum $M(\mu,\nu)$ on the Teichm\"uller space $\teich(\Sf)$.
\end{thm}

Our main result is the following, special cases of which we have already
proved in~\cite{KSQF} and~\cite{S1}:
\begin{thm}
\label{thm:mainiden} Let $\mu, \nu $ be 
two measured laminations which together fill up $\Sf$.
 Then as $\t \to 0$, the sequence
$G(\theta \mu, \theta \nu )$ of Theorem~\ref{thm:main} converges to
$M(\mu,\nu)$.
\end{thm}

 \noindent  The set of minima $M(\mu,t\nu)$ for $t \in (0,\infty)$ is a
line
$\L_{\mu,\nu} \subset \teich(\Sf)$ called the \emph{Kerckhoff line of
minima} of $\mu$ and $\nu$. 
Combining the above  results we obtain a complete proof of Conjecture
6.5 in~\cite{S1}:  
 \begin{thm}
\label{conj:mainconjecture} Let $\mu,\nu \in \ML$ be laminations which
fill 
up $\Sf$. Then 
the closure of $\P_{\mu,\nu}$ meets $\F$ precisely in $\L_{\mu,\nu}$.
 \end{thm}

 Notice that to prove the conjecture we need to invoke
Theorem~\ref{thm:otalbon} to construct the sequence 
$G(\theta \mu, \t \nu)$.
The Bonahon-Otal theorem  for rational laminations is based on the
Kerckhoff-Hodgson
theory of deformations of cone manifolds,  
 so our proof rests ultimately on the same thing. In~\cite{S1} we
were able to avoid this in some cases by directly proving  the existence
of a  sequence in $\P_{\mu,\nu}$ approximating a point $p \in 
\L_{\mu,\nu}$, however we required that the  supports of $\mu$ and
$\nu$ were pants decompositions and that 
a certain condition on the partial
derivatives of the lengths of these pants curves was satisfied at $p$.  
It would be nice to have a more general direct proof.

It is essential for the convergence in Theorem~\ref{thm:main} that the
bending measures $\pl^+$ and $\pl^-$ stay in bounded proportion. 
 For example, one might consider the case in which $\mu ,\nu $ are unit
measures $\delta_{\a}$ and $\delta_{\b}$ supported on fixed geodesics
$\alpha, \beta$
and study  the groups $G(\theta \delta_{\a},\phi \delta_{\b})$ with
$\phi/\theta \to 0$. In the case of a once punctured torus with $\alpha$
and $\beta$ a pair of generators,
one can check  by direct calculation (see Section~\ref{sec:opt}) that if
$\phi/\theta \to 0$ then this  sequence has no (algebraic) limit  as
$\theta \to 0$. 
We show that a similar phenomenon holds in general:  

\begin{thm}
\label{thm:converse} Let $\mu, \nu $ be 
two measured laminations which together fill up $\Sf$.
Then any sequence of groups $G(\theta \mu, \phi \nu)$ with $\theta, \phi
\to 0$ diverges (that is, no subsequence has an algebraic limit) unless 
$\theta/\phi$ is uniformly bounded away from $0$ and $\infty$. 
\end{thm}

 One also has to be careful if one wishes to allow 
$\mu, \nu$ to vary. It is easy to see by example that it is important
that the limit laminations themselves fill up $\Sf$.
We prove: 

\begin{thm}
\label{thm:diaglimit} Let $\mu, \nu $ be 
two measured laminations which together fill up $\Sf$, and suppose that
$\mu_n \to \mu$, $\nu_n \to \nu$ and $\theta_n \to 0$. Then the sequence
of groups
$G(\theta_n \mu_n, \theta_n
\nu_n)$ converges to $M(\mu,\nu)$ as $n \to \infty$.
\end{thm}

If the pleating loci $\pl^{\pm}$ are both rational, then, after removing
the pleating locus, one can double the convex core of the $3$-manifold
$\HH^3/G$ to obtain a cone manifold whose singular locus is the removed
bending lines. If the bending angle along an axis is $\phi$, then the
corresponding cone angle is $2(\pi- \phi)$. Thus one can regard the
above results as describing a special class of degeneration  of cone
manifolds to two-dimensional hyperbolic structures as all the cone
angles approach $2 \pi$ in a controlled way. 
Theorems~\ref{conj:mainconjecture} and~\ref{thm:converse}    give 
necessary and sufficient
conditions for such degeneration to occur.

\medskip

Of the above list, the  new results are
Theorems~\ref{thm:main},~\ref{thm:mainiden},~\ref{thm:converse} and~\ref{thm:diaglimit}.
The heart of the paper is Sections~\ref{sec:main}
and~\ref{sec:irrational},
in which we establish the following variant  of
Theorem~\ref{thm:main}:
\begin{prop}
\label{prop:mainrat} Let $\Sf$ be a hyperbolic surface of finite type,
and
suppose that  $\mu,\nu \in \ML(\Sf)$  
fill up $\Sf$. Then the groups $ G(\t\mu,\t\nu)$ lie in a relatively
compact set in $\QF$ and any accumulation point  as $\t \to 0$ is
Fuchsian. Moreover for any finite set  
$ \Gamma $ of simple curves on $\Sf$, there exists $c
>0,  $ such that  $|l_{\gamma } (p^+(G_{\t}))-
l_{\gamma } ( p^-(G_{\t}))| \le c   \t^2$ for all $\gamma \in  \Gamma$
and
all  sufficiently small $\t$.
\end{prop}

The paper is organized as follows. 
 After 
briefly summarising the background in 
Section~\ref{sec:preliminaries}, in Section~\ref{sec:identification} we
show that Theorems~\ref{thm:mainiden} and~\ref{thm:converse}  follow 
from Proposition~\ref{prop:mainrat}. One easily deduces
Theorem~\ref{thm:main} from  Proposition~\ref{prop:mainrat} and
Theorem~\ref{thm:mainiden}.  
In Section~\ref{sec:opt} we discuss the example of the
once-punctured torus referred to above.
In Section~\ref{sec:main} we prove Proposition~\ref{prop:mainrat} for
rational laminations and in Section~\ref{sec:irrational} in the general
case. Finally in Section~\ref{sec:diagonal}  we discuss diagonal limits,
showing  that there is sufficient uniformity in the estimates needed to
prove Proposition~\ref{prop:mainrat} to deduce
Theorem~\ref{thm:diaglimit}. 
  
 We should like to thank Vladimir Markovic and Young Eun Choi for encouragement, discussion and comments about the results in this paper.
   
 \section{Preliminaries} 
 \label{sec:preliminaries}
The background we need is mostly well known and explained at  
length elsewhere. Here we only give a brief summary and refer
to~\cite{S1} and elsewhere for more details.

Throughout the paper we write $g(\t) = O(\t)$ to mean that   $g(\t) \le
c\t$ for some fixed $c>0$ as $\t \to 0$.
We also write 
$g(\t) > O(\t)$ to mean that there exists $c>0$ such that 
$g(\t) > c\t$ as $\t \to 0$.

In general we shall be careful to specify the dependence of our
constants. Symbols $c,k$ and so on may denote different constants  in
different places but we label by a subscript if we need to refer back to
some particular earlier choice.

  \subsection{Quasifuchsian groups}
  \label{sec:quasifuchs}
  
  Let $\Sf$ be an oriented surface 
of negative Euler characteristic, homeomorphic
to a closed surface with at most a finite number of points removed. 
A {\em quasifuchsian group} $G$
   is the image of a discrete faithful representation 
  $\rho:\pi_1(\Sf)  \to PSL(2,\CC)$ 
  such that the   limit set of $G$ is a topological circle.
If $\Sf$ has punctures, we insist that   the images of loops around
boundary components are  parabolic.
 The limit set separates the regular set 
 into   two simply connected 
  $G$-invariant components $\Omega^{\pm}$ and each quotient  
  $\Omega^{\pm}/G$ is homeomorphic to $\Sf$.

 Two quasifuchsian groups are {\em equivalent} if the corresponding
representations are conjugate in $PSL(2,\CC)$.  {\em Quasifuchsian
space} $\QF(\Sf)$ is the space of equivalence classes.  It
has a   holomorphic structure induced from the
natural holomorphic structure of $SL(2,\CC)$.
A quasifuchsian group  is {\em Fuchsian} if the   limit set  is a round
circle. {\em Fuchsian
space} $\F= \F(\Sf)$ is the subset of $\QF(\Sf)$
corresponding to Fuchsian groups. 

As proved by Kerckhoff~\cite{KerckEA} p. 24, the induced real analytic 
structure on $\F$ is determined by the lengths of a finite number of
geodesics in $\pi_1(\Sf)$. These curves can always taken to be simple.
The arguments extend to show the same 
is true for the complex analytic structure on $\QF(\Sf)$,  if we replace
length by {\em complex length} $\lambda{(g)}$ given by the formula
$ \tr \rho(g) = 2 \ \ch \lambda{(g)/2}$ for $g \in \pi_1(\Sf)$.

For further details on these definitions, 
good references are~\cite{Marden, OtalH}.

\subsection{Geodesic laminations.}
\label{sec:laminations}
 
 Let $\Sf$ be a surface as above.
Given a hyperbolic structure on $\Sf$, a {\em geodesic
lamination} on
$\Sf$ is a closed union of pairwise disjoint simple complete geodesics
called its {\em leaves}.   
A {\em measured geodesic
lamination} $\mu$ consists of a geodesic
lamination, together with a transverse invariant measure on the leaves.
We denote the underlying lamination by $|\mu|$.
For reasons which will be clear in the next section, we only consider
 laminations with no leaves which end in a puncture.

We denote  the set of such
 measured laminations  by $\ML(\Sf)$.
This space 
is  topologised
 with the topology of weak convergence:
$ \mu_n \to \mu$ if $ \mu_n(T) \to \mu(T)$
for any transversal $T$. It is well known that $\ML(\Sf)$ is independent 
of the hyperbolic structure on $\Sf$.
For $\mu \in \ML$,  the length $l_{\mu}$ (relative to a given hyperbolic
structure on $\Sf$) is
   the total mass of the measure which
is the product of hyperbolic distance along the leaves of $|\mu|$ with
the transverse measure $\mu$.

 We call a measured geodesic lamination $\mu$ 
  {\em
rational} if its support $|\mu|$ consists entirely  of closed leaves,
and {\em irrational} otherwise.
(Note this is quite a different meaning from the term   {\em  
`arational'}   used for example in~\cite{OtalH}.) 
Let  $\S= \S(\Sf)$ 
denote  the set of all homotopy classes of simple closed 
non-boundary parallel curves on $\Sf$. If $\a_i$ are a set of disjoint
curves in $\S$, then 
by $\sum_i a_i\a_i $, 
 $a_i \in \RR^{+}$,  we mean  
 the measured lamination with support
$\cup_i \a_i$  which gives mass $a_i$ to each intersection
with  $\a_i$. 
Note that the maximum number of curves $\a_i$ in such a sum is
$3g-3+b$, where  $g$  is the genus of $\Sf$ and $b$ is the number of
punctures.
We   denote the
set of all rational measured laminations by
$\ML_Q(\Sf)$;
the set $\ML_Q$ is dense in $\ML$.

The length of the rational lamination $\sum_i a_i\a_i $
is just $\sum_i a_il_{\a_i} $, where $l_{\a_i}$ 
 is the hyperbolic length
of the geodesic $\a_i$.
Kerckhoff~\cite{KerckNRP, KerckEA} has shown that if
$\mu_n\in \ML_Q$ converges to $ \mu$ in $\ML$, then  
  $l_{\mu_n}$ converges to  $l_{\mu} $ 
 uniformly on compact subsets of $\F$, and hence is a real analytic
function on $\F$.
In a similar way, the geometric intersection number $i(\a,\a')$ of two
closed geodesics
$\a,\a'$ extends by linearity and continuity
to a continuous function
$i(\mu,\nu)$
on $\ML$, see for example~\cite{KerckNRP,OtalH}.

 Laminations $\mu,\nu \in \ML$ are said to {\em fill up} $\Sf$ if
 $i(\mu,\eta) + i(\nu,\eta) > 0$ for all $\eta \in \ML$.
An equivalent condition is that 
  every component of $\Sf - |\mu| \cup | \nu|$ 
  contains at most one puncture, whose closure, after 
  filling in the puncture if needed, is compact and 
  simply connected.

There is an obvious action of $\RR^+$ on $\ML$ given by 
 scalar multiplication $\mu \to t\mu$ for any $t >0$.
A {\em projective measured lamination} is an equivalence
class under this action. We write $[\mu]$ for the projective class of
$\mu$ and  denote the set of all  non-zero  projective measured
laminations by $\PML$. 
Thurston showed that $\PML(\Sf)$ can be viewed as the boundary of
$\F(\Sf)$: a sequence of structures $p_n \in \F$ converges to $  \xi \in
\PML$ if the lengths $\{l_{\g}(p_n)\}_{ \g \in \S}$ converge
projectively to the intersection numbers  $\{ i(\g,\xi) \}_{ \g \in
\S}$. It is not hard to deduce that if the laminations $\mu,\nu$ fill up 
$\Sf$ and if $p_n \in \F$ diverges, then at least one of $l_{\mu}(p_n)$
or $l_{\nu}(p_n)$ tends to $\infty$.

For more details on this material see for example~\cite{CEG}
or~\cite{OtalH}.

\subsection{The convex hull boundary and bending measures.}
\label{sec:convexhull}

 For any Kleinian group $G$,   let $\C = \C(G)$ be the hyperbolic  
convex hull  of the limit set  of $G$ in hyperbolic $3$-space $\HH^3$. 
If $G$ is quasifuchsian then 
$\bch$ has exactly two components $\bch^\pm$ which ``face'' the
components $\Omega^\pm$ of $\Omega$. The quotients
 $\bch^\pm/G$ are  homeomorphic
to $\Omega ^\pm/G$  and hence to $\Sf$.
(In the special case in which
$G$ is Fuchsian, ${\mathcal C}$ is
 contained in a single flat plane.
We regard  this as a degenerate case in which
 $\bch$ is two sided, each side facing one component
 of $\Omega(G)$.)

The structure of  $\bch$ is studied in detail in~\cite{EpM}.
Note that by convexity, $\bch$ must be embedded in $\HH^3$.
The ambient hyperbolic metric induces a metric on $\bch$ which endows
each 
component with its own hyperbolic metric; for a quasifuchsian group $G$
we shall denote the corresponding hyperbolic structures  in $\F(\Sf)$ by
$p^{\pm}(G)$. 
Each component of $\bch$ is the closure of a set of infinite sided ideal
polygons, each  contained in a hyperbolic plane  in $\HH^3$. These 
polygons, called the flat pieces of $\bch$,   are geodesic  not only  in
$\HH^3$ but also in the induced metrics on $\bch^{\pm}(G)$. (For a nice
picture of this, see \texttt{www.math.suny.edu/\~{}minsky}, reproduced
as Figure~12.6 in~\cite{Indra}.)
In each component of $\bch/G$, the closure of the complement of the flat
pieces is a geodesic lamination on $\Sf$, 
called the {\em bending lamination}, which 
 carries a   transverse
measure, the {\em bending measure}, denoted $\pl^\pm (G)$.

We note that no leaves of the   bending lamination
can limit on cusps of $\Sf$. For consider a   
horocycle of length $\epsilon$ round the cusp. 
The lift of the horocycle to $\HH^3$ bends by a definite amount
$\delta$  (fixed and independent of $\epsilon$) in every interval of
length $\epsilon$. 
By making $\epsilon$ sufficiently small, a comparison with the Euclidean
situtation shows that it is impossible for  $\bch$ to be embedded. 
This explains our assumption that laminations in $\ML$ 
contain no leaves which end in a puncture.

Each cusp on a hyperbolic surface is surrounded by a fixed area horocycle into
which no simple geodesic which does not end in the cusp penetrates. If
$\Sf_0$ is a hyperbolic structure on the surface $\Sf$,  let
$\Sf^C_0$ denote  the surface $\Sf_0$ minus
these fixed area horocyclic neighbourhoods of the punctures. 
The {\em non-cuspidal injectivity radius} $inj(\Sf^C_0)$ of $\Sf_0$ is
the 
radius of the largest embedded disc on $\Sf^C_0$.
By the above observation, a lower bound on $inj(\Sf^C_0)$ is equivalent
to a lower bound on the lengths of all simple curves on $\Sf_0$.

For a curve $\g \in \pi_1(\Sf)$ we always use
$\g^*$ and $\g^{\pm}$ to denote  the geodesic $\g$ in $\HH^3/G$ and its
geodesic representatives on $\bch^{\pm}/G$ respectively, and we
denote by $l_{\g^*}$ and $l_{\g^{\pm}} $   the  corresponding geodesic
lengths. Thus $l_{\g^*} \le l_{\g^{\pm}} $ and 
$l_{\g^*}= l_{\g^{+}} $ if $\g$ is contained in the support of $\pl^+$.
Notice that if $\mu \in \ML$
and $\pl^+(G) = \theta \mu$,  then the total  bending measure along
$\g^+$ is
$i(\gamma,\mu) \theta$. 

\medskip
 
We shall need the main result of~\cite{KSconvex}:
\begin{prop} \label{prop:KScont} The  maps 
$\QF \to \F, \ q \mapsto p^{\pm}(q)$
and  $\QF \to \ML, \ q \mapsto \pl^{\pm}(q)$
are continuous, where by definition $\pl^{\pm}(q) = 0$ if $q \in \F$. 
\end{prop} 
We remark that the results of~\cite{KSconvex} are stated  for
holomorphic families depending on one complex variable only, however the
theory of holomorphic
 motions extends to several variables~\cite{McS}
 and identical methods apply. 
 
    The following  straightforward result is \cite{S1}
Proposition 3.2, see also~\cite{BonO}.
 \begin{prop}
\label{prop:fillupbend} Let $G$ be quasifuchsian, $G \in \QF(\Sf)$.
Then the bending measures  $\pl^{\pm}(G)$ 
 fill up $\Sf$.
\end{prop}
We remark that the proof in~\cite{S1} is not quite complete in the case
in which $\pl^{\pm}$ are irrational, because we  omitted the possiblity
that  a complementary region of 
the union of the two laminations is simply connected or a once punctured
disk but non-compact.
However in this case one obtains a semi-infinite geodesic $\a$ in the
complement of  both $|\pl^{\pm}|$.  The accumulation points of $\a$ form
a geodesic lamination contained in both 
$\bch^{+}$ and $\bch^{-}$. Following the same idea as in~\cite{S1},
this is easily seen to be impossible. For more details,
see~\cite{KerckLM} Lemma 4.4.

\paragraph{Rational bending laminations.}
The structure of a component of $\bch$ is particularly simple when its
 bending measure is rational, 
say $\pl^{+} = \sum_i  \t a_i \a_i$ for some $\t>0$. In this case  $\bch^+$ consists of
pieces of hyperbolic planes
which meet exactly along the lifts of  axes of the curves $ \a_i$,
in such a way that the exterior angle of intersection 
(i.e. the angle {\em outside $\C$}) along an axis which projects to $
\a_i$ is $\t a_i$. Since $\bch $ is convex,
 all the angles
$\t a_i$ have the same sign.
They are measured so that $\t = 0$ exactly when  the oriented planes
containing the adjacent flat pieces coincide. In particular,    $G$ is
Fuchsian if and only if $\pl^{+} = 0$ (so that $\pl^{-} = 0$ follows
automatically).

We remark that if the bending lamination is rational, then each flat
piece of $\bch$ faces a disk in the regular set which contains a Cantor
set of limit  points in its boundary. Such `ghost circles' are a highly 
visible feature of many limit set pictures, see for
example~\cite{Indra}.

 \subsection{Earthquakes and Quakebends.}
\label{sec:bending}

 The time $t$  left earthquake  along 
a  lamination $\mu \in \ML$ is a real analytic map 
 $\E_{\mu}(t): \F \to \F$ which generalises  the classical
Fenchel-Nielsen twist. Let $T$ be a transversal to $|\mu|$
with endpoints in distinct complementary components of the
 $|\mu|$.  The earthquake shifts the component on the right a distance
$t\mu (T)$ relative to the one on the left, inducing a new
hyperbolic  metric $\E_{\mu}(t)(p) $ on an initial hyperbolic structure
$p
\in \F$. In particular, if $\mu = \sum_i a_i \a_i$,  then for each $i$,
the earthquake $\E_{\mu}(t)$
twists
by hyperbolic distance $ta_i \a_i$ around the closed geodesic $\a_i$. 

The map $(p,t) \mapsto \E_{\mu}(t)(p)$ is a flow on  $\F$
 which induces a tangent vector field 
$\p_{t_{\mu}}$.
In \cite{KerckNRP}, Kerckhoff showed that if $\nu \in \ML$, then the
length
$l_{\nu}$ is a
 real analytic function of $t$
along the flow, strictly convex if $i(\mu,\nu)>0$ 
and constant otherwise. Wolpert~\cite{Wolp} proved the 
famous antisymmetry relations
$\dd l_{\nu}/\dd t_{\mu} = -  \dd l_{\mu}/\dd t_{\nu}$.

Complexifying the parameters corresponds to passing from 
Fuchsian to quasifuchsian groups. In this context 
the earthquake $\E_{\mu}(t)$ has a natural extension
 to a {\em left quakebend}  $\E_{\mu}(\tau), \tau \in \CC$.
We shall only need the construction  relative to a hyperbolic
structure $p_0$ corresponding to an initial Fuchsian group $G_0$, in
which form  it is explained in detail in~\cite{EpM}. 
In addition to  shifting complementary components of $|\mu|$ on $p_0$
through a relative distance  $\Re \tau \mu (T)$
as above, the map 
 $\E_{\mu}(\tau)$  bends
the righthand component  through the angle $\Im \tau \mu (T)$
  relative to the left one.
This deforms the group $G_0$, given by a representation
$\rho_0: \pi_1(\Sf) \to  PSL(2,\RR)$, into a group 
$\E_{\mu}(\tau)(G_0)$ given by a representation $\rho(\tau): \pi_1(\Sf)
\to  PSL(2,\CC)$. A quakebend    with purely imaginary parameter $\tau
\in i\RR$
is called a {\em pure bend}.

The following result is~\cite{KSbend} Theorem 8.8.
In~\cite{KSbend} it is explained in the context of a punctured
torus but the proof clearly extends to a general surface.
\begin{prop}
\label{prop:smallbends} 
Let $G_0 \in \F$ and  $\tau \in \CC$. Then provided $|\Im \tau|$ is
sufficiently small, $\E_{\mu}(\tau)(G_0)$ is
quasifuchsian. If 
$\E_{\mu}(\tau)(G_0)$ is quasifuchsian, then
$p^+(\E_{\mu}(\tau)(G_0)) = p^+(G_{0})$
and $\pl^+(\E_{\mu}(\tau)(G_0)) = (\Im \tau)\mu$.
\end{prop}

 By Lemma 3.8.1 of~\cite{EpM}, the groups
$G_{\tau \mu}$  
depends  holomorphically on $\tau$.
One can therefore define the complex derivatives
$\dd l_{\nu}/\dd {\tau}_{\mu}$; 
complex versions of the antisymmetry formulas have been proved by
Kourouniotis~\cite{Kour}.

 \section{Identification of the limit}
 \label{sec:identification}
   
 In this section we show that  that Theorems~\ref{thm:mainiden}
and~\ref{thm:converse}  follow  from Proposition~\ref{prop:mainrat}. 
The proof is based on Kerckhoff's result which we stated as
Theorem~\ref{thm:kerck}. In fact  he proved a rather stronger
statement:
if  $\mu$ and $\nu$  fill up $\Sf$, then 
$ l_{\mu} + l_{\nu}$  has a unique {\em critical point} on $\F$.

  As usual, let 
$G(\t) = G(\t \mu, \t \nu)$ denote a quasifuchsian group for which
$\pl^+ = \t \mu$ and  $\pl^- = \t \nu$, and let $p^{\pm}(\t) $ denote
the Fuchsian structures on $\bch^{\pm}/G(\t)$.

\begin{prop} \label{prop:iden} Let $\mu, \nu $ be 
two measured laminations which together fill up $\Sf$.
Suppose that as $\t_n \to 0$, the groups $ G(\theta_n \mu, \theta_n \nu
)$ converge  to  $p \in \F$ and in addition 
that for any finite set   $ \Gamma  \subset \S$, there exists $c
>0$, depending on $\Gamma$,  such that  $|l_{\gamma } (p^+(\t)-
l_{\gamma } ( p^-(\t)| \le c   \t^2$ for all $\gamma \in  \Gamma$ and
all  sufficiently small $\t$.
Then $p=M(\mu,\nu)$.
\end{prop}
(We remark that the bound $|l_{\gamma } (p^+(\t)-
l_{\gamma } ( p^-(\t)| \le c   \t$ would be sufficient but our work
leads naturally to $\t^2$.) 

\noindent \begin{proof}
By the above noted extension to
Theorem~\ref{thm:kerck}, it is sufficient to show that $p$
is a critical point of the length function $l_{\mu}+l_{\nu}$ on $\F$.
Since by~\cite{KerckLM} Theorem 3.5 the tangent vectors $\p_{t_{\xi}},
\xi \in \ML$ span the tangent
space $T_p(\F)$ to $\F$ at $p$, it is enough to show that 
$$\frac{\dd l_{\mu}}{\dd t_{\xi}}(p) + \frac{\dd l_{\nu}}{\dd
t_{\xi}}(p) =0  \ \  \forall \xi \in \ML,$$
and hence, by the antisymmetry of the derivatives, that 
$$\frac{\dd l_{\xi}}{\dd t_{\mu}}(p) + \frac{\dd l_{\xi}}{\dd
t_{\nu}}(p) =0  \ \ \forall \xi \in \ML.$$

The proof of~\cite{KerckEA} Lemma 2.4  shows that
$l_{\xi}$ is a  real analytic function on $\F(\Sf)$. 
 By a straightforward extension of the arguments, see also~\cite{KSQF}
Theorem 6.3, one sees that  
$l_{\xi}$ extends to  a complex analytic function on $\QF$.
As noted above, for fixed $p_0 \in \F$ and $\mu \in \ML$, 
the group 
$\E_{\mu}(\tau)(p_0)$ obtained by  quakebending by $\tau$ along $\mu$
depends
holomorphically on $\tau$; in particular
$l_{\xi}(\E_{\mu}(\tau)(p_0))$ is an analytic function of  $\tau$.

Now a quasifuchsian group is completely
determined by the Fuchsian structure $p^+$ on $\bch^{+}$ and the bending
measure $\pl^+$. This leads to the key observation that we can reach
$G(\t
\mu, \t \nu)$ either by starting at 
$p^{+}(\t) $ and making the pure bend  $\E_{\mu}(i \t)$ through the
angle $\t$ along $\mu$, or by 
starting at $p^-(\t)$ making a pure bend through   $-\t$ along $\nu$. 
The idea is to use this to make Taylor series expansions of the length
of an arbitrary geodesic in $G(\theta \mu, \theta \nu ) $ in two ways
and  compare
the results.

Write $q(\t) $ for $G(\t \mu, \t \nu)$ and let $\sigma^+: [0,1]\to \QF$
be the pure bend path between $p^+(\t) =
\sigma^+(0)$
and $q(\t) =  \sigma^+(1)$, so that $\sigma^+(t) = \E_{\mu}(it \t)$ is
the
quasifuchsian group obtained by bending $p^+(\t)$
through the angle $t\t$ along $  \mu$.  Expanding from 
$p^{+}(\t) $ we obtain:

\begin{equation}
\label{eq:topexpand}
l_{\xi}(q(\t)) = l_{\xi}(p^+(\t)) + i \t \frac{\dd l_{\xi}}{\dd
t_{\mu}}(p^+(\t)) - \t^2 \frac{\dd^2 l_{\xi}}{\dd t_{\mu}^2}(r^+(\t))
\end{equation}
where $r^+(\t) \in \sigma^+$. 
With a similar definition of $\sigma^-$, expanding from $p^{-}(\t) $ we
get:
 \begin{equation}
 \label{eq:botexpand}
 l_{\xi}(q(\t)) = l_{\xi}(p^-(\t)) - i \t \frac{\dd l_{\xi}}{\dd
t_{\nu}}(p^-(\t)) - \t^2 \frac{\dd^2 l_{\xi}}{\dd t_{\nu}^2}(r^-(\t))
\end{equation}
where $r^-(\t) \in \sigma^-$.

Now by Proposition~\ref{prop:KScont}, $\lim_{n \to \infty} p^{\pm}(\t_n)
=p$,
so the points 
$r^{\pm}(\t_n)$ lie in some compact neighbourhood of $p$ in $\QF$. It
follows that the second derivatives in the above expressions are
uniformly
bounded as $n \to \infty$. Equating imaginary parts we find 
\begin{equation}
\label{eq:bottop}
\frac{\dd l_{\xi}}{\dd t_{\mu}}(p^+(\t_n)) + \frac{\dd l_{\xi}}{\dd
t_{\nu}}(p^-(\t_n)) = O(\t_n) 
\end{equation}
with constants depending on $\xi$.

Now choose the set $\Gamma$ to be a finite set of curves which determine
the analytic structure on $\F$.
Since $l_{\xi}$ is real analytic on $\F$,
our hypothesis gives 
\begin{equation}
\label{eq:bottopnu}
\frac{\dd l_{\xi}}{\dd t_{\nu}}(p^+(\t)) -\frac{\dd l_{\xi}}{\dd
t_{\nu}}(p^-(\t)) = O(\t^2). 
\end{equation}

It follows that 
\begin{equation}
\frac{\dd l_{\xi}}{\dd t_{\mu}}(p^+(\t_n)) + \frac{\dd l_{\xi}}{\dd
t_{\nu}}(p^+(\t_n)) = O(\t_n).
\end{equation}
Using Proposition~\ref{prop:KScont} again,  we make the
required
conclusion
by taking limits as $n \to \infty$.
\end{proof}

 \noindent {\sc Proof of Theorem~\ref{thm:mainiden}.}
 This follows immediately from Proposition~\ref{prop:iden} and
Proposition~\ref{prop:mainrat} (to be proved in Sections~\ref{sec:main}
and~\ref{sec:irrational}).

\begin{remark} {\rm Theorem 5.1 in~\cite{S1} is effectively the special
case
of Proposition~\ref{prop:iden} in which $\mu$ and $\nu$ are both
rational laminations
whose supports are pants decompositions of $\Sf$. However the proof
in~\cite{S1} fails completely in the irrational case. Notice also that
the above proof rests heavily on the  assumption  that the limit
group $p$ exists; even if $p^{\pm}(\t)$ remain close the error terms
$\frac{\dd^2 l_{\xi}}{\dd t_{\mu}^2}(r^+(\t)),\frac{\dd^2 l_{\xi}}{\dd
t_{\nu}^2}(r^-(\t)) $ might well become unbounded if $p^{\pm}(\t) \to
\dd
\F$. (Wolpert's formula~\cite{Wolp} for the second
derivatives shows that these terms contain a factor $1/l_{\xi}$.)} 
 \end{remark}

 \medskip

 Theorem~\ref{thm:converse} is proved by a similar method. This time we 
only need 
Proposition~\ref{prop:propA} from Sections~\ref{sec:main}
and~\ref{sec:irrational}.

\medskip

 \noindent {\sc Proof of Theorem~\ref{thm:converse}.}
Suppose as usual that  $\mu, \nu \in \ML$   fill up $\Sf$, and
suppose that we have a sequence of groups
$G_n = G(\theta_n  \mu, \phi_n  \nu)$ with $\theta_n, \phi_n
\to 0$ for which $\theta_n/\phi_n \to 0$.
 We have to show that  the sequence $G_n$ diverges. If not, then
(by passing to a subsequence if necessary) we may suppose that $G_n$ has
an algebraic limit  $ G_{\infty}$.

We claim  that    $G_{\infty}$ is Fuchsian.
To see this, write $\pl^{\pm}_n, \pl^{\pm}_{\infty}$ for the bending
measures of $ G_n,  G_{\infty} $ respectively.
By~Proposition~\ref{prop:KScont},   $\pl^{\pm}_n \to
\pl^{\pm}_{\infty}$. 
Since  our  hypothesis implies that  $\pl^{\pm}_n   \to 0$,  we deduce
that 
$\pl^{\pm}_{\infty} = 0$. Thus each of  $\bch^{\pm}(G_{\infty})$ is
contained in a single hyperbolic plane, from which it follows that the
regular set of $G_{\infty}$ contains at least two circular invariant
domains. By Accola's theorem a group with at least two simply connected
invariant domains is quasifuchsian; if these domains are circular 
it must be Fuchsian.

Write $p^{\pm}_n, p^{\pm}_{\infty}$ for the Fuchsian structures on
$\bch^{\pm}(G_n)/G_n, \bch^{\pm}(G_{\infty})/G_{\infty}$ respectively.
Following exactly the method used to arrive at
equation~(\ref{eq:bottop}) in   Proposition~\ref{prop:iden},  
we find  \begin{equation}
\label{eq:bottop1}
\t_n \frac{\dd l_{\xi}}{\dd t_{\mu}}(p^+_n) + \phi_n \frac{\dd
l_{\xi}}{\dd
t_{\nu}}(p^-_n) = O(\t_n^2 + \phi_n^2).
\end{equation}

As before, let $\Gamma$   be a finite set of curves which determine
the analytic structure on $\F$. 
In a compact neighbourhood of $G_{\infty}$, the non-cuspidal injectivity
radii of 
the structures $p^{\pm}_n$ are uniformly bounded below and the lengths
of the curves in $\Gamma$ are uniformly  bounded above. Thus
Proposition~\ref{prop:propA} gives that
$|l_{\g^*}(G_n) - l_{\g^+}(G_n)| \le O(\t_n^2)$ and 
 $|l_{\g^*}(G_n) - l_{\g^-}(G_n)| \le O(\phi_n^2)$.
(Here as elsewhere, $l_{\g^*}$ and $l_{\g^{\pm}} $ denote respectively
the lengths of the geodesic $\g$ in $\HH^3/G$, and of its geodesic
representatives on $p^{\pm}$.)
Hence, noting that $ l_{\g}(p^+_n)  = l_{\g^+}(G_n)$, we have 
$|l_{\g}(p^+_n) - l_{\g}(p^-_n)| \le O(\t_n^2+ \phi_n^2)$ for all
$\g \in \Gamma$. Combining this with the fact that $l_{\xi}$ is a real
analytic function on $\F$,
it follows that  for any $\xi \in \ML$,
\begin{equation}
\frac{\dd l_{\xi}}{\dd t_{\nu}}(p^+_n)  + \frac{\dd l_{\xi}}{\dd
t_{\nu}}(p^-_n)  = O(\t_n^2+ \phi_n^2). 
\end{equation}

Together  with~(\ref{eq:bottop1}) this gives
\begin{equation}
\t_n \frac{\dd l_{\xi}}{\dd t_{\mu}}(p^+_n) + \phi_n \frac{\dd
l_{\xi}}{\dd
t_{\nu}}(p^+_n) = O(\t_n^2+ \phi_n^2).
\end{equation}

Dividing through by $\phi_n$ and taking limits we deduce 
(again using Proposition~\ref{prop:KScont}) that 
$$
 \frac{\dd l_{\xi}}{\dd
t_{\nu}}(p_{\infty}) = 0.
$$

Since this holds for all $\xi \in \ML$ we deduce that 
$l_{\nu}$ has a critical point at $p_{\infty} \in \F$ which
(since not all derivatives  along all possible earthquake paths
$\E_{\xi}(t)$ can vanish, see e.g.~\cite{KerckLM} p.194) is impossible.
\qed

\section{A special example} 
 \label{sec:opt} 
 Before proceeding to the proof of Proposition~\ref{prop:mainrat}, we
pause to
examine one of the few examples  in
which one can write down exact formulae for the relationship
between bending angles and lengths and hence explore the limit behaviour
explicitly.\footnote{For some other examples in which explicit
caluations can be made, see~\cite{DSExamples}.} Namely, take $\Sf$ to be
a once-punctured torus and let
$\alpha, \beta \in \pi_1(\Sf)$ intersect exactly once. Thus
$\pi_1(\Sf)$ is the free group generated by $\alpha$ and $\beta$ and 
the commutator $\alpha \beta \alpha^{-1} \beta^{-1}$ represents a loop
around the puncture.   
We shall study  the case in which $\pl^+  \in [\alpha]$ and $\pl^- \in 
[\beta]$, in other words, the surfaces $\bch^{\pm}$ are bent
 along  axes which project to 
$\alpha$ and $ \beta$ respectively.  Although this situation appears to
be very special, quite similar geometry appears in the general case.

By~\cite{KStop}  Lemma 4.6, 
if a geodesic  $\g$ is  contained in the  bending lamination  then its
image  in $PSL(2,\CC)$   has real trace. 
As shown in~\cite{PS}, this  gives the equations
\begin{eqnarray}
\label{eq:angleislength}
\cos \t_{\a}/2 = \ch l_{\b^*}/2 \; \th l_{\a^*}/2, \ \ 
\cos \t_{\b}/2 = \ch l_{\a^*}/2 \; \th l_{\b^*}/2 
\end{eqnarray}
 relating the bending angles   $\theta_{\a}, \theta_{\b}$  to the 
 lengths $l_{\a^*}, l_{\b^*}$ of  the geodesic representatives  $\a^*$
and $\b^*$ of  $\a,\b$ in $\HH^3/G$. 
Moreover suitably chosen  lifts  $\ax A,
\ax B$ of  $\a^*$ and $\b^*$ are mutually
perpendicular at distance $d$, where
\begin{equation}
\label{eq:optdist}
\ch  d \; \sh l_{\a^*}/2 \; \sh  l_{\b^*}/2 = 1.
\end{equation}
(Here  $A,B \in G$ are translation  by $l_{\a^*}, l_{\b^*}$ along the
respective axes.)

\begin{figure}[hbt]
\begin{center}
\caption{Configuration of bending axes on the once punctured torus.}
\label{quad}
\end{center}
\end{figure}

Since $\pl^+ \in [\a]$, the axis $\ax A$ lies on $\bch^+$.
  Let $\tilde \b^+$ be the lift of $ \b^+$ whose endpoints on $\partial
\HH^3$ are the same as those of $\ax B$, 
 and let  
 $P = \ax A \cap \tilde \b^+$. Let $Q$ be the foot of the 
perpendicular from
$P$ to  $\tilde \b^*$, see Figure~\ref{quad}. 
 Since
$B$ is purely hyperbolic, the quadrilateral $\cal Q$
  with vertices $P$, $Q$ and $B(P), B(Q)$ is planar. Note that $|PQ|=d$, 
and that  since 
 $\Pi$ is orthogonal to $\ax A$, the line $P B (P)$ makes an angle
$(\pi-
\t_{\a})/2$ with $PQ$.
By applying the 
  quadrilateral formulae (see~\cite{Beardon}  Theorem 7.17.1) to $\cal
Q$ we obtain:
 \begin{equation}
 \label{eq:plane1}
 \sh {d} = \cth l_{\b^*}/2 \; \tan \t_{\a}/2,   \ \ \ 
 \ch d \; \sh l_{\b^*}/2= \sh l_{\b^+}/2.
 \end{equation}

\begin{lemma} Suppose that   $a,b>0$  are fixed and that  $d= d(\t)$ is
the distance between
    $\ax A$ and $\ax B$  in the group $G(a\t,b\t)$. 
Then   $d \le
O(\t)$ as $\t \to 0$.
\end{lemma} 
\begin{proof}  Our hypothesis means that $ \t_{\a} = a \t$ and $\t_{\b}
= b \t$. Suppose first that there is some subsequence along which 
$l_{\b^+} \ge c> 0$. 
If in addition  $l_{\b^*} \ge c'> 0$, then equation~(\ref{eq:plane1}) 
gives $d \le O(\t)$.

Otherwise, passing to a further subsequence, we may assume that 
$l_{\b^*} \to  0$.
From~(\ref{eq:plane1})  we have  
$$ \th d = \frac{\tan \t_{\a}/2}{\th l_{\b^*}/2} \; \frac{\sh
l_{\b^*}/2}{\sh l_{\b^+}/2},$$
from which, since $l_{\b^+}$ is bounded away from $ 0$, it follows that
$ d \le O(\t)$.

 By interchanging the roles of $\a$ and $\b$, we conclude that either
both  $l_{\b^+} \to 0$ and   $l_{\a^-} \to
0$; or  
$d \le O(\t)$ as $\t \to 0$.  Suppose that the first alternative
applies. Then certainly also $l_{\a^*} \to 0$. However $l_{\a^*}$ and
$l_{\b^+}$ are the geodesic lengths of $\a$ and $\b$ on the Fuchsian
structure $\bch^+/G$, and by the collar lemma this situation is
impossible. 
 \end{proof}

\begin{cor} Let $\Sf$ be a once punctured torus with generators $\a$,
$\b$, and suppose that  $\mu = a\delta_{\a}, \nu =b\delta_{\b}$.  Then
 the group $G(\t \mu,\t\nu)$ converges to a Fuchsian group 
as $\t \to 0$. Moreover
the hypotheses of Theorem~\ref{thm:mainiden} hold, so that  the limit is
the minimum on $\F$ of the function $a l_{\a} + b l_{\b}$.  
\end{cor}
\begin{proof} From the lemma we have that
$d \le O(\t)$ as $\t \to 0$.  We deduce from~(\ref{eq:plane1}) 
and its analogue with $\a$ and $\b$ interchanged   that both $l_{\a^*}$
and $l_{\b^*}$ are bounded 
away from $0$, and then from~(\ref{eq:optdist}) that they are both
bounded above.  This is sufficient to ensure (up to a subsequence)  the
existence of the algebraic limit of $ G(\t \mu,\t \nu)$. (One way to
see this is to use the Markov equation which relates $\tr AB $ to $\tr
A$ and $\tr B$.) One also sees that $l_{\a^*}/l_{\a^+} \to 1$ and
similarly for $\b$. Moreover not only the axes $l_{\a^*}$ and
$l_{\b^*}$, but also $l_{\a^*}$ and $l_{\b^+}$, and $l_{\a^-}$ and
$l_{\b^*}$ are orthogonal. This is enough to ensure that the limit of
each of the two Fuchsian structures $\bch^{\pm}/G$ also exists and
equals
the limit  of $G(\t \mu,\t \nu)$. The remaining details are left to the
reader.   
 \end{proof}

In the above discussion we made crucial use of the fact that 
 $\t_{\a}/\t_{\b}$ is bounded away from $0$ and $\infty$  (in fact
constant) as $\t \to 0$.  Without this hypothesis, the result fails.
In fact rearranging equation~(\ref{eq:angleislength})  one obtains 
\begin{equation}
 \label{eq:lengths}
\sh l_{\a^*}/2 = \sin \t_{\b}/2 \; \cot \t_{\a}/2,  \ \ 
\sh l_{\b^*}/2 = \sin \t_{\a}/2 \; \cot \t_{\b}/2.  
\end{equation}
  If only one of $\t_{\a}$ and $\t_{\b} $ converges to $ 0$ then
 one of $l_{\a^*}$ and $l_{\b^*}$ diverges to $\infty$. 
If both $\t_{\a}$ and $\t_{\b}  $ converge  to $ 0$
then 
 $\sh l_{\a^*} \sim \t_{\a}/\t_{\b}$ and $\sh l_{\b^*} \sim
\t_{\b}/\t_{\a}$. If the ratio is unbounded either above or below,
again at least one of $l_{\a^*}$ and $l_{\b^*}$ diverges to $\infty$.
Note,
however, that  in this case
$1/ \ch d =   \sh   l_{\a^*}\sh l_{\b^*} \to 1$ so we still  get that 
$d \to 0$.

\section{The main limit theorem} 
 \label{sec:main}
 In this section we establish Proposition~\ref{prop:mainrat} in the case
in
which $\mu$ and $\nu $ are rational. 
The idea of the proof is as follows. First, in
Proposition~\ref{prop:propA}, we
establish an upper bound $d \le O(\t)$ for the
distance between any point  on the lift of a closed geodesic to
$\bch^{\pm}$ and the corresponding axis in $\HH^3$, under the hypothesis
that the length of the corresponding curve on $\bch^{\pm}/G$
is bounded below. Our estimate also controls the ratio of the lengths
on $\bch^{\pm}/G$ and in $\HH^3/G$. 
Then in Proposition~\ref{prop:lowerbound} 
\label{page:stages} we prove a \emph{lower} bound $d \ge O(\t)$ for the  
distance between any point
on a bending line and the opposite side of $\bch$.
In Proposition~\ref{prop:lengthupperbound} we play off these two bounds
against each other to deduce an upper bound on the lengths of all
bending lines. This is sufficient to establish the existence of the
limit. Another use of Proposition~\ref{prop:propA} also 
establishes the necessary estimate on the variation of length of curves
in $\Gamma$.

\begin{prop}
\label{prop:propA} Fix  $L>0$. Let  
 $\mu \in \ML$ be fixed
 and suppose that $G \in \QF(\Sf)$ is such  that $\pl^+(G) = \theta
\mu$.  
For any  $\gamma \in \pi_1(\Sf)$, let   $\tilde \gamma^+$, $\tilde
\gamma^*$ be lifts of $  \gamma^+$, $ 
\gamma^*$ to $\HH^3$ with the same endpoints on
$\partial \HH^3$.
Suppose that $l_{ \gamma^+} \ge L$. Then there is a 
universal constant $
\t_0$, and a constant
 $c_0=c_0(L)$, such that 
 for all $\theta < \t_0$ and any $P \in \tilde  \gamma^+$:
$$d( P, \tilde \gamma^*) \le c_0 i(\gamma,\mu) \theta \ \ {\rm and } \ \
l_{ \gamma^*} \ge (1-c_0(i(\gamma,\mu) \theta)^2) l_{ \gamma^+}.$$
\end{prop}

In this section we prove this result on the assumption that $\mu $
is rational. The extension to the general case is not hard and is done
at the beginning of Section~\ref{sec:irrational}. 

\smallskip
\begin{remark}   
{\rm    
\begin{enumerate} 
\item   This result certainly has applications beyond the present one.
It
will be
noted in the proof that the Kleinian group $G$ can be quite general and
that all that is
needed is that $\tilde \g^+$ lie on a pleated surface with bending angle
$O(\t)$; convexity is also not required.
\item It is crucial in our statement that the constants $c_0(L)$ and $
\t_0$
do not depend on the hyperbolic structure of $\bch^+$. If $\mu$ is
multiplied by a scalar $t>0$, then the term  $i(\gamma,\mu)$ scales
accordingly. 
 \item  The result fails without the hypothesis of a lower bound on
$l_{\gamma^+}$. In this sitation, the distance 
between $l_{ \gamma^+}$ and $l_{ \gamma^*}$ may become infinite with no
control on the ratio $l_{ \gamma^+}/ l_{ \gamma^*}$. This   can be seen
by letting $\t_{\a} \to 0$ and $\t_{\b} \to \pi/2$ in
equation~(\ref{eq:lengths}) in Section~\ref{sec:opt} and then 
examining~(\ref{eq:plane1}). 
\item The following variant  has been
proved independently by Lecuire~\cite{Le} (without the estimate
of the distance between $l_{\gamma^+}$ and $l_{\gamma^*}$): \\ 
{\it Suppose $\epsilon \le \pi/12$ and $i(\g,\mu) \le \epsilon$. Then 
$ l_{\g^+} \le (1+\tan \epsilon) (l_{\g^*}+6\epsilon)$.} 
\end{enumerate}
} 
\end{remark}

 \begin{figure}[hbt]
\begin{center}
\caption{Configuration for Proposition~\ref{prop:propA}.}
\label{fig2}
\end{center}
\end{figure}

Let  $\g \in \pi_1(\Sf)$ have lifts $\tilde \g^+$, $\tilde \g^*$ as in
the statement of the proposition.  Pick $P \in \tilde \g^+$
and let $\hat   \g$ denote the
piecewise geodesic arc in $\HH^3$
joining the  points $\g^n(P), n \in \ZZ$,  see Figure~\ref{fig2}.
(By abuse of notation we are also using $\g$ to denote the element of
$G$ in the conjugacy class of $\g$ which fixes the axis $\tilde \g^*$.)
The idea is to  estimate the distance between $\tilde \g^+$
and $\hat \g$, and then between $\hat \g $
and $ \tilde \g^*$, at the same time comparing  their lengths.
To do this we use two lemmas about piecewise geodesic arcs in $\HH^3$
based on a nice idea in~\cite{CEG}
Theorem 4.2.12.  First, a simple result about hyperbolic triangles.

\begin{lemma}
\label{lemma:angleofp} Let $ABC$ be a hyperbolic trangle with exterior
angle $\phi$ at $C$. If $h = d(C, AB)$ then $\th h \le  \sin \phi$.
\end{lemma}
\begin{proof} Let $X$ be the foot of the perpendicular from $C$ to $AB$.
Since both of the angles $ACX$ and $BCX$ are less than $\pi/2$, the line 
through $C$ and perpendicular to $CX$ is outside the triangle $ABC$ and
hence makes an angle $ \psi < \phi$ with $CB$.

Let $E$ be the point at which the extension of  $CB$ meets   $\partial
\HH^2$, and let $Z$ be the foot of the perpendicular from $E$ to the
extension of $CX$. The extension $\l$ of $AB$ divides $\HH^2$ into two
half planes; since $EZ$  cannot cut $\l$, it must lie on the side not
containing $C$ so that $|CZ| \ge |CX|$.
By the angle of parallellism formula,  $ \th |CZ| = \sin \psi$. 
The result follows.
  \end{proof}

Now for the estimates  based on~\cite{CEG}.\footnote{We note that the
last sentence in the statement of~\cite{CEG}
Theorem 4.2.12 is incorrect. The proof however  is correct and our
version here indicates one way of proving what was clearly intended.}
The following notation is convenient. 
Let $\sigma$ be   {\em any}  piecewise geodesic arc   in $\HH^3$
with endpoints $X$ and $X'$. 
For $P \in \sigma$, let $v(P) = v(P,\s)$ be the (positive) angle at $P$
between the forward vector
along $\s$ at $P$ and the forward vector along the line extending  
$XP$ pointing away from $X$. 
Likewise let $w(P)= w(P,\s )$ be the angle at $P$ between the backwards
vector along $\s $ at $P$ and the forward vector along the line
extending   $X'P$ pointing away from $X'$. 
The configuration is shown in Figure~\ref{fig3}.

\begin{figure}[hbt]
\begin{center}
\caption{Configuration for Lemma~\ref{lemma:ceg}.}
\label{fig3}
\end{center}
\end{figure}
 
\begin{lemma} 
\label{lemma:ceg} Let $\sigma$ be   {\em any}  piecewise geodesic arc  
in $\HH^3$
with endpoints $X$ and $X'$, and let $\hat \s$ be the $\HH^3$ geodesic
joining $X$ to $X'$.  Suppose that for all $P \in \s$, 
both the angles 
 $v(P,\s)$, $w(P,\s)$  are bounded above by $\phi$. Then 
$ l_{\hat \s} \ge (\cos \phi ) l_{\s} $ and $\th d(P, \hat \s) \le \sin
2\phi$ for   all $P \in \s$, where $l_{\s} $ and $ l_{\hat \s}$ are the
lengths of $\s$ and $\hat \s$ respectively. 
\end{lemma}
\begin{proof}  
Suppose the arc $\s$ has successive 
 bends at points $X=X_0$, $ X_1$, $X_2, \ldots, X_k=X'$.
Let us
compare the geodesic distance $x=|XP|$ with the distance
$t=  \\ \sum_{i=0}^r |X_iX_{i+1}| + | X_rP|$ measured along the broken
arc
 $\s$. Obviously, $x$ is a piecewise $C^1$ function of $x$, only
failing to be differentiable at the bends $X_i$. It is not hard  to
check that on each open arc,  $\frac{dx }{dt} = 
\cos v(P,\s)$.
Thus the first part of the  result follows by integrating (using
the obvious continuity of $x$ as $P$ moves 
through  each point $X_i$).

For the second part, our hypothesis implies that the exterior angle at
$P$ of the triangle $PXX'$ is at most $2 \phi$.  An application of
Lemma~\ref{lemma:angleofp} gives the result.
\end{proof}

\begin{lemma}
\label{lemma:ceg2} Fix $L, k>0$.  Suppose $  \s$ is a piecewise geodesic
arc in $\HH^3$ for which each segment has length at least $L$ and the
angle at each bend is in absolute value  at most $k \t$.
Then there exist   $\t_1>0$, and $c_1=c_1(L)$ depending only on $L$,  
such that
for
any $P \in  \s$, and all
 $\theta < \t_1$, the angles $v(P,  \s)$ and $w(P,  \s)$ 
are at most $c_1k \theta$. 
\end{lemma}
\begin{proof} This can be deduced as in~\cite{CEG} Theorem 4.2.12; for
convenience we give slightly different version here.
We work with $v(P)$; the argument for $w(P)$ is the same.
As before, suppose $\s$ meets successive 
 bending lines in points $X=X_0$, $ X_1$, $X_2, \ldots, X_k=X'$.
Let  $P $ be a point on the open segment  $X_rX_{r+1}$, so that 
$v(P)=v(P,\s)$ is the (positive) angle from the oriented line 
$X_rX_{r+1}$ and the extension of $XP$ through $P$, oriented away from
$X$. Clearly $v(P)$ decreases as $P$ moves along the arc from 
$X_r$ to $X_{r+1}$. Let $u(X_{r+1})$ denote the angle at $X_{r+1}$
between the forward vector along the extension of $X_rX_{r+1}$ and
$XX_{r+1}$,
so that $v(X_{r+1})$ is the sum (or difference if the bend is negative)
of  $u(X_{r+1})$
and the bending angle at $X_{r+1}$.  Writing  $t = d(X_r,P) $   and  $v
= v(P)$, then  as shown in in~\cite{CEG}, 
 $d v /dt = - \sin v / \th {|X P|}$.  
 Hence  we can estimate the decrease in $v$ along $X_rX_{r+1}$ by 
$$v(X_r) - u(X_{r+1}) = - \int_{X_{r+1}}^{X_r} \frac{d v}{dt} \  dt \ge
L\sin
u(X_{r+1})   \ge Lu(X_{r+1}) /2$$ provided the initial angle  $v(X_r)$
is chosen less than some fixed $\t_1$ for which $\sin \t_1 \ge \t_1/2$
say.
Thus  $$u(X_{r+1}) \le  2v(X_r)/(2+L) $$ and so,  by our hypothesis on
the bending angles,  $$v(X_{r+1}) \le  2v(X_r)/(2+L) + k \t.$$

Choose $c_1(L)=(2+L)/L$. Suppose inductively that the  angle 
$v(P)$ is at most $c_1k \t$ for any point $P$ in the closed subsegment
of  
$\s $ between   $X$ and $X_r$.  
This is certainly true for $r=0$ since on the first segment $v(P) = 0$.
By hypothesis $v(X_r) \le k \t$ so as above, $v(X_{r+1}) \le  2c_1k
\t/(2+L) + k \t  \le c_1 k\t $
by our choice of $c_1$. The result follows. 
\end{proof}

  \noindent {\sc Proof of Proposition~\ref{prop:propA}.}
Assume that $\mu$ is rational. As illustrated in Figure~\ref{fig2}, pick
a point $P \in \tilde \g^+$ and denote one full translation length of
$\tilde \g^+$ along
$\bch^+$   from  $P$ to $\g(P)$  by $\s^+$.
(Here $\g$ denotes the particular choice of element
in the conjugacy class of $\g$ which fixes $\tilde \g^+$.) We  first
compare $\s^+$ to the $\HH^3$ geodesic  $\hat \s$
joining $P$ to $\g(P)$.
Since $\s^+$ is
geodesic on $\bch^+$, the angle between successive geodesic segments of
$\s^+$ is
bounded
above by the bending angle on $\bch^+$ (see for example~\cite{KSbend}
Lemma 6.2).  Let $Q$ be a
point on one such segment  $X_rX_{r+1}$. By  Gauss-Bonnet  the angle 
between
$PQ$ and  $X_rX_{r+1}$ is bounded above by $\t i(\g,\mu)$, as is the
angle
between 
$\g(P) Q$ and $X_rX_{r+1}$. 
By Lemma~\ref{lemma:ceg}, 
$$ l_{\hat \s} \ge (1-2(i(\mu,\g) \t)^2 ) l_{\s^+} \ \ {\rm and} \ \ 
d(Q, \hat \s) \le 2i(\mu,\g)\t,$$   where $ l_{\hat \s}$ is the length
of $ \hat \s$, that is, the  distance in $\HH^3$ between $P$ and
$\g(P)$.

Now let $\hat \g_n$ denote the piecewise geodesic arc $\cup_{r=-n}^{n-1}
\g^r(\hat \s)$, so that $\hat \g_n$ joins the points $  \g^r(P), 
\g^{r+1}(P)$, for $r = -n,
  \ldots, n-1$. 
Applying Lemma~\ref{lemma:ceg2}  to 
$\hat \g_n$, we see that there exists $c_1>0$, depending only on
$L$,  such that for any $Q \in \hat \g_n$, and all
sufficiently small $\theta$, the angles  between $\g^{-n}(P)Q$ and
$\hat \g_n$, and between $Q\g^{n}(P)$ and $\hat \g_n$, are at most
$c_1i(\mu,\g) 
\theta$.  

Applying Lemma~\ref{lemma:ceg} again shows that any point on the
geodesic arc joining 
$\g^{-n}(P)$ to $\g^{n}(P)$ is within distance $2c_1i(\mu,\g) \t$ of
$\hat \g_n$
and that 
$$d_{\HH^3}(\g^{-n}(P),\g^{n}(P)) \ge  2n( 1 - 2(c_1i(\mu,\g)  \t)^2 )  
l_{\hat \s}.$$

Since  $\g^{-n}(P)$ and  $\g^{n}(P)$ converge to the negative and
positive fixed points of the axis $\tilde \g^*$ respectively, the arc
joining 
$\g^{-n}(P)$ to  $\g^{n}(P)$ converges  to $\tilde \g^*$. The result
follows.
 \qed

 We now turn to the lower bound from p.~\pageref{page:stages}.
A {\em support plane} $\Sigma$ to $\bch$ is a complete
hyperbolic plane in $\HH^3$ which meets $\C$ and such that all of $\C$
is contained in one of the two half spaces cut out by $\Sigma$. 
 We repeatedly use the following easy fact. 

\begin{lemma}
\label{lemma:dontmeet} Let $\Sigma^+$, $\Sigma^-$ be support planes to 
$\bch^+$, $\bch^-$ respectively. Then $\Sigma^+ $ and $\Sigma^-$ are
disjoint.
\end{lemma}
\begin{proof} Denote by $D^{\pm}$ the open disks in $\Chat$ which are
the ends at infinity of the half spaces  in $\HH^3$ cut out  by
$\Sigma^{\pm}$ and not containing $\C$.
 We claim that neither disk $D^{\pm}$  contains any limit points of $G$.
In fact if there were a limit point  in $D^{+}$,  then by
convexity there would be a line segment  contained in $\C$  and meeting
the half planes on both sides of $\Sigma^{+}$, which is impossible.

Since  $G$ is quasifuchsian, its regular set has two
components $\Omega^{\pm}$. By the definition of $\bch^{\pm}$, we have 
$D^{\pm} \subset \Omega^{\pm}$. 
If $\Sigma^+$ and $ \Sigma^-$ meet, then so do $D^+$ and $ D^- $,
contradicting the fact that $\Omega^+ $ and $\Omega^-$ are disjoint.
\end{proof}
 
\begin{figure}[hbt]
\begin{center}
\caption{Support planes to $\bch^+$ and $\bch^-$ illustrating
Proposition~\ref{prop:lowerbound}.}
\label{fig4}
\end{center}
\end{figure}

The idea is that there is not only a maximum but also a {\em minimum}
distance between a bending line on one side of $\bch$  and any support
plane of the other side, because, as illustrated in Figure~\ref{fig4}, 
a pair of support planes which meet at
a definite angle are at a definite distance away from any other plane
which is disjoint from them both. 

\begin{prop}
\label{prop:lowerbound} There exist universal constants  $c_2>0, \t_2>0$
with the property that if
 the point $P$ is on a bending line $\tilde \a$ in $\bch^+$, then $d(P,
\bch^-) \ge c_2 i(\mu,\tilde \a)\t$ for all  $\t < \t_2$. 
\end{prop}   
 \begin{proof} We need to estimate the nearest possible approach to
$\bch^+$ of a support plane to $\bch^-$. Write  $k$  for the transverse
$\mu$-measure
$i(\mu,\tilde \a)$ of $\tilde \a$, that is, the weight of
the axis $\tilde \a$ in the lamination $\mu = \sum_i a_i \a_i$.  Let
$\Sigma, \Sigma'$ be the support planes which meet along $\tilde \a$.
First
consider the situation in the plane $\Pi$ through $P$ and orthogonal to
$\tilde \a$. Let $\lambda, \lambda'$ be the lines in which  $\Sigma,
\Sigma'$
meet $\Pi$, and let $\zeta$ be the half-line starting at $P$ which
bisects the
angle $\pi - k\t$ between $\Sigma$ and $\Sigma'$.  We claim there
are no points of $\bch^-$ on $\zeta$ within distance $kO(\t)$ of $P$. 

If this is not true, then there is a support plane of $\bch^-$ which
meets $\zeta$. By
Lemma~\ref{lemma:dontmeet}, any such plane to $\bch^-$ is disjoint from
$\Sigma$ and $ \Sigma'$. The limiting case is that in which such a  
plane meets $\Pi$ in the line $\eta$ which joins the ends of 
$\lambda, \lambda'$ at infinity, forming a triangle with two ideal
vertices and exterior angle $ k\t/2$ at $P$.
If $h$ is the perpendicular distance from $P$ to 
$\eta$, then by the angle of parallelism formula, $ \th h = \sin ( 
k\t/2)$. The claim follows.

 Now we look in the plane $\Pi'$ orthogonal to $\Pi$, which contains the
bending axis $\tilde \a$ and the half-line $\zeta$.  Let $Q$ denote the
intersection
of $\zeta$ with $\eta$, so that $|PQ| = h$. 
Let   $\Delta$ be  the triangle with vertices the two ends of $\tilde
\a$ at
$\infty$ and $Q$.
We claim that no support plane of $\bch^-$ meets $\Delta$. In fact any
such support plane meets $\Pi'$ in a line $\lambda''$; by the
first
part of the proof $\lambda''$ does not meet the segment $PQ$, nor, by
Lemma~\ref{lemma:dontmeet}, does it meet $\tilde \a$. Therefore if 
$\lambda''
\cap
\Delta \neq \emptyset$, $\lambda''$ must enter and exit $\Delta$ across
the same side, which is impossible. 

Now we calculate the radius of the maximal half disk centre $P$ 
contained in  $\Delta$. Let $\pi - 2\phi$ be the interior angle at $Q$. 
By the angle of parallelism formula, $\sin \phi = \tanh h$ so that by
the
above, $ \phi = k\t/2$. The required radius is the perpendicular
distance  $h'$ from $P$ to either of  the two other sides of $\Delta$,
and hence $\sh h' = \cos \phi \; \sh h$, from which we deduce that $h'
\ge kO(\t)$.

Finally we consider the intermediate case of a plane $\Pi''$ containing
the line $\zeta$ and making an angle between $0$ and $\pi/2$  to the
bending
axis $\tilde \a$. 
We consider the quadrilateral $\Q$ with two sides $\Sigma \cap \Pi''$,
$\Sigma' \cap \Pi''$ which meet at $P$, and whose remaining two sides
are
the lines through $Q$ meeting $\Sigma \cap \Pi''$ and $\Sigma' \cap
\Pi''$ 
on $\partial \HH^3$. Just as above, we argue that no support plane of
$\bch^-$
can intersect the interior of $\Q$. Then we calculate the maximal radius
$h''$ of a disk centre $P$ whose intersection with the half planes cut
off by $\Sigma \cap \Pi''$, $\Sigma' \cap \Pi''$ and containing $\zeta$
is
contained in $\Q$. 

Let $\pi - 2\psi, \pi - 2\phi''$ be the angles in $\Q$ at $P$, $Q$
respectively.  As before, $\sh h'' = \cos \phi'' \; \sh h$  so it is
enough to show that $\phi'' \le k O(\t)$. Let $Q'$ be the foot of the 
perpendicular from the one vertex of $\Q$ on $\partial \HH^3$ to $PQ$.
Then $\sin
\phi'' = \th |Q'Q|  \le \th h = kO(\t)$ which
proves the result.
\end{proof}

We can now play Proposition~\ref{prop:propA} and
Proposition~\ref{prop:lowerbound} against each other to get a bound on
the lengths on the bending lines.  The idea is that if 
a very long segment of the geodesic representative $\tilde \a^-$ of a
bending line $\a \subset |\mu|$ is entirely contained in a flat piece of
$\bch^-$, then it must get  very close to the actual bending line
$\tilde \a^* = \tilde\a^+$ on $\bch^+$, contradicting
Proposition~\ref{prop:lowerbound}.
We need the following lemma about skew quadrilaterals which is 
 Theorem 2.4.6 in~\cite{CEG}. It follows from an explicit
calculation of the distance between points on two fixed geodesics, and
the fact that the distance function is convex. 

\begin{lemma}
\label{lemma:cegdistance}   Let $X_1X_2Y_2Y_1$ be a skew hyperbolic
quadrilateral. 
Suppose that $d(X_i,Y_i) \le v$ for $i=1,2$. Let $\eta >0$ be
given and let $Z$ be any point on   $X_1X_2$ with  $d(X_i,Z) \ge \eta$
for each $i$. Let $u = d(Z,Y_1Y_2)$ be
the distance from $Z$ to the line $Y_2Y_1$. Then $\sh u \le \sh
v/\ch \eta$. 
\end{lemma}

\begin{prop}
 \label{prop:lengthupperbound}
 The lengths  $l_{\mu^{\pm}}, l_{\nu^{\pm}}$ of the bending laminations
on $\bch^{\pm}/G(\t)$ are uniformly bounded above as $\t \to 0$.
 \end{prop}
 \begin{proof} 
It will be  enough to show that there is a uniform upper bound on 
$l_{\alpha^-}$ for any component  $\a$ of
$|\mu|$. For by similar reasoning we also obtain an upper bound on
$l_{\beta^+}$ for any component  $\b$ of
$|\nu|$, and hence {\em a fortiori} on $l_{\beta^-}$.

A lift $\tilde \alpha^-$ of $\alpha^-$ on $\bch^-$ is  partitioned into
a finite number  of geodesic segments by the points at which it meets
the bending lines $|\nu|$ on $\bch^-$. We claim that
the length  of each segment is
uniformly bounded above as $\t \to 0$. As usual, let $\tilde
\alpha^{\pm}$ be lifts of $ \alpha^{\pm}$ with the same endpoints on
$\partial \HH^3$.

We may as well suppose that $l_{\alpha^-} \ge 1$, so  we can
apply Proposition~\ref{prop:propA} to $\bch^-$ and $\alpha$, with $L =
1$, to show that $d(P, \tilde \alpha^+) \le O(\t)$ for all $P \in \tilde
\alpha^-$.
Let $X_1,X_2$ be successive points at which  $\alpha^-$
meets $|\nu|$, and let $Y_1,Y_2$ be the feet of the
perpendiculars from $X_1$ and $X_2$ to  
$ \tilde \a^+$.  We shall apply Lemma~\ref{lemma:cegdistance} to the
skew
quadrilateral $X_1,X_2, Y_2,Y_1$. 
We have just shown that  $d(X_i,Y_i) \le O(\t)$.  
Let $Z$ be the midpoint of $Y_1Y_2$ and let $u = d(Z,X_1X_2)$. 
Since $\tilde  \a^*$ is a bending line and since the segment from $X_1$
to $
X_2$ is contained in $\bch^-$,  Proposition~\ref{prop:lowerbound} 
gives   $ u \ge O(\t)$.   
Applying Lemma~\ref{lemma:cegdistance}, we find
$\sh {u} \le O(\t) / \sh {y}$ where  $y = |Y_1 Y_2|/2$. Given the lower
bound on $u$, this  is
impossible if $y \to \infty$. We deduce that $y$ is uniformly bounded
above. Since $Y_1Y_2$ is the perpendicular projection of $X_1X_2$ 
 through a distance   $O(\t)$,  we also get a bound on $X_1X_2$.
Summing over all segments  gives a uniform upper bound on $l_{\a^-}$ as
required. 
 \end{proof}

 \begin{cor} 
 \label{cor:compactset} The structures $p^{\pm}(\t)$ lie in a compact
set 
in $\teich(\Sf)$.
 \end{cor}
 \begin{proof} This is~\cite{ThuII} Corollary 2.3. If the
structures $p^+(\t)$ were not in a
compact set in $\teich$, then we could find a subsequence
converging to a point $\xi$ in the Thurston boundary  $\PML$. Since the
systems $|\mu|$ and $|\nu|$ together fill up the surface, $\xi$ has
non-zero intersection
number with at least one component  $\delta$ of either $ |\mu|$ or
$|\nu|$, and hence 
$l_{\delta^+} \to \infty$ as $\t \to 0$. This contradicts the above
proposition and proves the
claim.
 \end{proof} 
We can now  prove the main result of this section.

\medskip
 \noindent{\sc Proof of Proposition~\ref{prop:mainrat}.}
By Corollary~\ref{cor:compactset}, for small $\t$ all 
of the structures $p^{\pm}(\t)$ lie in a compact set $K$ in
$\F$. 
Choose a sequence  $\t_n \to 0$ along which  
$p^{+}(\t_n) \to p^+_{\infty}$ and $p^{-}(\t_n) \to p^-_{\infty}$  for 
points $p^{\pm}_{\infty} \in K$. 

By compactness,  there is a uniform
lower bound to the non-cuspidal injectivity radius of all surfaces in
$K$, equivalently, a uniform lower bound to length of all simple
geodesics.
Therefore we may apply
Proposition~\ref{prop:propA} to see that   $l_{\g^*} \ge
(1-O(\t_n^2))l_{\g^{\pm}}$  for   any 
curve $\gamma \in \S$, where the constant  depends only
on $i(\g,\mu) $, $ i(\g,\nu)$ and $K$.  Writing $G_n$ for $G(\t_n)$, we
deduce that 
$ 1-O(\t_n^2) \le l_{\g^-}(G_n)/ l_{\g^+}(G_n) \le  1+O(\t_n^2) $ as $n
\to
\infty$. Since in  $K$  the lengths $l_{\g^{\pm}}(G_n)$ are also
uniformly
bounded above,  we deduce that $|l_{\g^+} (G_n) - l_{\g^-} (G_n)| \le
O(\t_n^2)$ as $ n \to \infty$, with constant depending only on $\g$.

Applying this to  a fixed finite collection of curves $\g_i$ whose
lengths determine the complex analytic structure on $\QF$, we deduce in
particular that $p^+_{\infty}= p^-_{\infty}$. We also deduce  
also from the expansions in equations~(\ref{eq:topexpand})
and~(\ref{eq:botexpand}) in Proposition~\ref{prop:iden}  that $\lim_{ n
\to \infty} l_{\g_i} ( G_n) =
l_{\g_i}(p^+_{\infty})$ and hence that the algebraic limit of the groups
$G_n$ is
the Fuchsian group corresponding to the structure $p^+_{\infty}$.
(Since  $p^+_{\infty}$ is certainly geometrically finite, the limit is
strong.)\qed

Combining this result with Proposition~\ref{prop:iden},  we see   that
the limit in Proposition~\ref{prop:mainrat} is
in
fact the minimum $M(\mu,\nu)$ and hence independent of the subsequence
chosen. This completes the proof of Theorems~\ref{thm:main}
and~\ref{thm:mainiden}  in the rational case.

\smallskip
   Notice that the compact set $K$ in the above proof is {\em not} given
uniformly in terms of $\t$, but depends in an unspecified way on  the
point $p^+_{\infty}$. This will cause us some grief in
Section~\ref{sec:diagonal}.

  \section{Extension to the irrational case} 
 \label{sec:irrational}
 We now discuss the extension of
Proposition~\ref{prop:mainrat} to the  case in which $\mu,\nu$ are
irrational. The  proof in the last section does not
immediately  extend mainly because the constants involved in the final
estimates  in Proposition~\ref{prop:lengthupperbound}   depend heavily
on the number of bending lines in $|\mu|$ and $|\nu|$.
In order  to control the `size' of  a measured lamination $\xi \in \ML$,
 fix once and for all a set $\Gamma = \{ \g_1,\ldots, \g_k \} \subset
\S$ of
curves
which fill up
$\Sf$, and  set $||\xi||_{\Gamma} = \sum_j i(\g_j,\xi)$. Since the
curves fill up, we have $||\xi||_{\Gamma} >0$.  Notice that 
$||\xi||_{\Gamma} >0$ is independent of the hyperbolic structure on
$\Sf$.

\medskip 
We begin by completing the proof of Proposition~\ref{prop:propA} for
irrational $\mu$. 
\medskip

\noindent {\sc Proof of Proposition~\ref{prop:propA}.} All we need to do
is adapt the first part of the proof
 from Section~\ref{sec:main} to the case in which $\mu \notin
\ML_{\QQ}$.  
As before, pick a point $P \in \tilde \g^+$ and denote the arc  of
$\tilde \g^+$ along
$\bch^+$   from  $P$ to $\g(P)$  by $\sigma^+$.
 We  want to compare $\sigma^+$ to the $\HH^3$ geodesic  $\hat \sigma$
joining $P$ to $\g(P)$.

Recall that the bending measure and distance along any arc $\kappa$   on
$\bch^+$ is defined in terms of  finite approximations called  `roofs'.
Namely we approximate $\kappa$ by the geodesic segments which join  the
intersection points of $\kappa$ with a finite number of leaves of $\mu$,
and
 measure the arc length and total bending angle along this finite
approximation in the obvious way.  The distance and bending angle along
$\kappa$ are by definition the infima, over all possible roofs, of the
corresponding finite approximations, see~\cite{EpM} and
also~\cite{KSconvex}. 

 In the current situation, we obtain the required type of estimate for
any roof exactly as before and  the required comparisons 
$$ l_{\hat \s} \ge (1-O((i(\mu,\g) \t)^2) ) l_{\s^+} \ \ {\rm and} \ \ 
d( P, \hat \s) \le O(i(\mu,\g) \t) $$  for any $P \in \s^+$  follow.
The remainder of the proof is exactly as before. 
 \qed

 We need the following extension of Proposition~\ref{prop:lowerbound}.
Although we keep the same names, the constants involved are not exactly
the same as those in the earlier version.
 
\begin{prop}
\label{prop:lowerbound1}  There exist $\epsilon_2, c_2,\t_2>0$ with the
following
property:   if  $P \in \bch^+$ lies on a geodesic segment $\s$ of length
at most $\epsilon_2$,  then $d(P,\bch^-) \ge c_2i(\s,\mu)\t$ for all $\t
< \t_2$. 
\end{prop}

The idea of this result is clear but a careful proof requires some work. 
The following two lemmas
 control  changes of angle as we move short distances along a  piecewise
geodesic arc.

\begin{lemma}
\label{lem:idealangledecrease1} Let $Z \in \partial \HH^3$ and let $\l$
be an oriented line with endpoints distinct from $Z$.  For $i=1,2$ let  
$X_i$ be points on $\l$ and let $\phi_i$ be the positive angle at $X_i$
between the forward direction of $\l$ and the extension of $ZX_i$
through
$X_i$.   Then $\phi_2/\phi_1 > 1- |X_1X_2|$.  
\end{lemma}
\begin{proof} Let $W$ be the foot of the perpendicular from $Z$ to $\l$.
For any point $X \in \l$, define $\phi$ as in the statement and let 
$t = |WX|$. By the angle of parallelism formula, $\th t = \cos \phi$.
Differentiating we find $ \frac{d \phi} {dt} = - \sin \phi$. (Note this
is the limiting case of a similar formula used in the proof of
Lemma~\ref{lemma:ceg2}.) Integrating along the arc from $X_2$ to $X_1$
gives
$ \phi_1 - \phi_2 \le  |X_1X_2| \sin \phi_1$, from which the result
follows.  
\end{proof}

Now let $\s$ be a piecewise geodesic arc in $\HH^3$ with a finite number
of bends $X_0, \ldots, X_k$ and endpoints $Z,Z'$  in  $\partial \HH^3$.
As in the discussion just before Lemma~\ref{lemma:ceg}, 
for $P \in \sigma$, let $v(P) = v(P,\s)$ be the (positive) angle at $P$
between the forward vector
along $\s$ at $P$ and the forward vector along the line extending  
$ZP$ pointing away from $Z$.

\begin{lemma}
\label{lem:idealangledecrease2} Suppose that $\s$ is a piecewise
geodesic arc on $\bch^+$ with initial and final points $Z,Z' \in
\partial \HH^3$ and successive bends at points $  X_0, X_1,\ldots, X_k
\in \HH^3$. Suppose the angle between successive segments at $X_i$ is
$\phi_i$. 
Then with the notation above, 
$$v(X_r) \ge (1- \sum_{i=0}^{r} |X_iX_{i+1}|) (\sum_{i=0}^{r} \phi_i).$$  
\end{lemma}
\begin{proof}  We prove this by induction on $r$.
For $r=0$, the result follows from the definitions. 
Assume $r>0$ and that the result holds for $r-1$.  Setting $\e_i =
|X_{i-1}X_i|$, this
means that 
$$v(X_{r-1}) \ge \Bigl(1- \sum_{i=0}^{r-1} \e_i \Bigr)
\Bigl(\sum_{i=0}^{r-1} \phi_i \Bigr).$$

Let $\psi_i$ denote the angle between the extension of $ZX_i$ and
$X_{i-1} X_i$, so that $v(X_i) = \psi_i + \phi_i$. 
Using Lemma~\ref{lem:idealangledecrease1} with $\l$ the geodesic
extending  the arc $X_{r-1}X_r$, we find $\psi_{r}/v(X_{r-1})  \ge (1-
\e_r) $. Hence
 $$v(X_{r}) =  \phi_r + \psi_r \ge   \phi_r  + (1- \e_r)\Bigl(1-
\sum_{i=0}^{r-1}
\e_i \Bigr)\Bigl(\sum_{i=0}^{r-1} \phi_i\Bigr).$$
This last expression is easily seen to be greater than
$ (1- \sum_{i=0}^{r} \e_i )(\sum_{i=0}^{r} \phi_i)$ as required.
\end{proof} 

Finally we need a lemma  to control the angles at which the arc $\s$
intersects bending lines.

\begin{lemma}
\label{lem:disjtleaves} Choose $\epsilon < \cosh^{-1} \sqrt 2$. Suppose
that $\l$,
$\l'$ are disjoint lines
in the hyperbolic plane $\HH$  (possibly  meeting on $\partial \HH$), 
and that $P \in \l, P' \in \l'$ are such that $|PP'| \le \epsilon$. 
Then the line through $P$ orthogonal to $\l$ meets $\l'$ in  a point
$Q$;
moreover
  $|PQ| < O(\epsilon)$  and $\angle PQP' \ge \pi/2 - O(\epsilon)$.
\end{lemma}
\begin{proof} If the orthogonal to $\l$ through $P$ does not meet $\l'$
then $d(\l,\l') \ge \cosh^{-1} (\sqrt 2)$, this number being  the
altitude of
a triangle with angles $\pi/2,0,0$. This proves the first statement.
  
Write $|PQ|=x$ and  $\angle PQP' = \phi$. Let $\phi_0$ be the angle
between $PQ$ and the line joining $Q$ to the endpoint of $\l$ on
$\partial \HH$ on the
same side of $PQ$ as $P'$.  Clearly $\phi \ge \phi_0$ and by the angle
of
parallelism formula, $ \sin \phi_0 = 1/ \ch x$.  

Let $P''$ be the foot of the perpendicular from $P$ to $\l'$; then $h =
|PP''| \le \epsilon$.
By trigonometry in triangle $PQP''$ we have $\sin \phi = \sh h / \sh x$.
Combining these observations we find $\sh h \ge \th x$,  from which it
follows that $x \le O(\epsilon)$ and hence that $\phi \ge \pi/2 -
O(\epsilon)$ as claimed.
\end{proof}

\noindent{\sc Proof of Proposition~\ref{prop:lowerbound1}} 
First suppose that $\mu$ is rational. Let $\s$ be a goedesic segment in
$\bch^+$ of length $\e$ to be determined later.
We begin by showing that we may assume that $\s$ is more or less
orthogonal
to all the bending lines. Let the first and last bending lines cut by
$\s$ be $\l', \l''$ respectively. 
 
We claim we can always find an $\HH^3$ geodesic  $\l$ through $P$ and
completely contained in $\bch^+$. This is obvious if  $P$ is on a
bending line. If not,  there is
some $\HH^3$ geodesic  $\l$ through $P$ contained in a flat piece of
$\bch^+$ and disjoint from all leaves of $|\mu|$ (possibly
meeting leaves of $|\mu|$ on $\partial \HH^3$).

Let $\Pi$ be the plane through $P$ orthogonal to $\l$.   
We shall first show that we may replace $\s$ by the segment $\s_1$ in 
$\Pi \cap \bch^+$ joining $\l'$ to $ \l''$.
In fact  applying  Lemma~\ref{lem:disjtleaves}  we see that  $\s_1$
meets the same
leaves as $\s$  (so that $i(\s_1,\mu) = i(\s,\mu)$) and that the length
of $\s_1$ is at most $O(\epsilon)$.
Thus we may as well work entirely in the plane $\Pi$.

As usual, let 
$X=X_0, X_1, \ldots, X_k= X'$ denote the  points at which $\s_1$
meets the bending lines of $\bch^+$. 
Let $\phi_i$ denote the angle at $X_i$ between the segments 
$X_{i-1}X_i$ and $X_i X_{i+1}$, and let $\t_i$ be the angle between the
support planes which meet at $X_i$. We claim that we may as well replace
$i(\s,\mu) = \sum \t_i$ by 
$ \sum \phi_i$.
In fact from Lemma~\ref{lem:disjtleaves} we see that $ \s_1$
is almost orthogonal to the bending line through $X_i$, crossing at the
angle $\psi_i$ say. These angles are related by the formula
$\tan  \phi_i/2 = \tan \t_i/2   \sin \psi_i$.
We deduce $\phi_i
> (1-O(\epsilon)) \t_i$.  
By definition, $\sum_{i=0}^{k} \t_i = i(\s, \t \mu) $.  
Therefore $\sum_{i=0}^{k} \phi_i > i(\s,\mu)\t/2 $, say, for all
sufficiently small $\epsilon$.

Now let $X_{-1}$ and $X_{k+1}$
respectively be the bending points  immediately  preceeding $X_0$ and
immediately 
following $X_k$ on the extension of the $\bch^+$ geodesic containing
$\s_1$, and let $Z,Z'$ be the points where the continuations
of $X_0X_{-1}$ and $X_kX_{k+1}$ meet $\partial \HH^3$. Suppose that $P$
is on the arc $X_{r-1}X_{r}$, $0< r\le k$. As above, if $P$ is not on a
bending line we may insert an extra line $\l$ containing $P$ and
disjoint
from the other lines in $|\mu|$ and treat $\l$ as a bending line with
bending angle $0$. Applying   
Lemma~\ref{lem:idealangledecrease2}, we get $v(P) > (1-\epsilon)
\sum_{i=0}^{r-1} \phi_i$ and similarly $w(P) > (1-\epsilon)
\sum_{i=r}^{k}
\phi_i$, where $w(P)$ is the angle   at $P$  between the forward vector
along the line extending   $Z'P$ pointing away from $Z'$
and the backwards direction along $\s_1$. Thus the angle $v(P) + w(P)$
between the lines $ZP$ and $PZ'$
at $P$ is at least $(1-\epsilon) i(\s,\mu) \t/2$.

We can now complete the proof more or less exactly as in
Proposition~\ref{prop:lowerbound}, replacing the support planes which
meet along the bending line by the planes containing the lines $ZP$ and
$PZ'$  and orthogonal to $\Pi$. Examination of
Proposition~\ref{prop:lowerbound} shows that we only need check that no
support plane $\Sigma^-$ to $\bch^-$ meets either of these lines in
$\HH^3$. (In
the plane orthogonal to $\Pi$, the line $\l$ is already in  $\bch^+$
and can be treated as a bending line.) If $\Sigma^-$ meets $PZ$ in
$\HH^3$, then
$Z$  is contained in the open disk on $\partial \HH^3$ spanned by
$\Sigma^-$ and containing no limit points (see
Lemma~\ref{lemma:dontmeet}).  But $Z$ is also on the boundary of any
support plane to $\bch^+$ which contains the arc  $X_{-1} X_{0}$.
By Lemma~\ref{lemma:dontmeet}, this is impossible.

Finally, we need to deal with the case in which $\mu$ is irrational.
Since none of the above estimates depend on the number of bending lines
which meet $\s$, approximating $\mu$ by finite laminations as explained
in the part of proof of Proposition~\ref{prop:propA} at the beginning of
this section will work. 
\qed

As in the rational case,  we are  now set to play the upper and lower
bounds against each other. We need the following lemma on intersection
numbers. 
\begin{lemma}
\label{lemma:intbound0} Suppose that  $\mu$ and $\nu$ fill up $\Sf$.
Then there exists
$c_3>0$ such that $i(\g ,\mu) +i(\g ,\nu) > c_3$ for all $\g \in \S$.
\end{lemma}
\begin{proof} Since the result depends only on intersection numbers we
can work entirely   with a fixed hyperbolic structure $p_0 \in \F$ whose
non-cuspidal injectivity radius is $\rho_0>0$ say. 
If the result is false, then we can find a subsequence $\gamma_n \in \S$
such
that $i(\g_n,\mu) +i(\g_n,\nu) \to 0$. Passing to a further subsequence
we may assume that there is a sequence $h_n >0 $ such that $h_n
\delta_{\g_n}  \to  \xi $ in $\ML$ where $l_{\xi}(p_0) =1$.  Then
$h_nl_{\g_n} \to 1$ and since $l_{\g_n} \ge \rho$  we have $h_n
\le 2/\rho_0$. Thus  
$i(h_n\g_n,\mu) +i(h_n\g_n,\nu) \to 0$ and hence, taking limits,
$i(\xi,\mu) +i(\xi,\nu) = 0$. Since $\mu$ and $\nu$ fill up $\Sf$, this
is impossible. 
\end{proof}

Now we can establish a uniform lower bound to the non-cuspidal
injectivity radii of the structures $p^{\pm}(\t)$.

\begin{prop}
\label{prop:injradius} Suppose that $\mu,\nu \in \ML$ and that $G(\t)=
G(\t \mu, \t
\nu) \in \QF(\Sf)$, and let $p^+(\t)$ denote  the Fuchsian structure on
$\bch^{+}/G(\t)$.  Then there exist  $\rho_*>0$ and $ \t_3>0$ such that
$l_{\g}(p^+(\t)) > \rho_* $  for all 
 $\g \in \S$ and all $\t < \t_3$.
\end{prop}
\begin{proof} 
Use the Margulis lemma to choose  $\rho>0$ such that
that if $\delta_1,\delta_2 \in \S$ and $i(\delta_1,\delta_2)>0$, 
then $l_{\delta_i} > \rho>0$ for at least one $i$.
Further reducing  $\rho$ if necessary, we may also assume that
$\rho < \epsilon_2$,   chosen as in
Proposition~\ref{prop:lowerbound1}.

Suppose that $\omega \in \S$ is such that $l_{\omega}(p^+(\t)) < \rho$
for some $\t$.  Since the curves in $\Gamma$ fill up $\Sf$,  we must
have $i(\omega, \delta) > 0$
for some  $\delta = \delta(\t) \in \Gamma$ for which 
$l_{\delta}(p^+(\t)) > \rho$.
Applying Proposition~\ref{prop:propA} to $\delta$ and $\bch^+(\t)$, we
see that $d(P, \tilde \delta^*) \le c_0(\rho)i(\delta,\mu)\theta$ for
all $P \in \tilde 
\delta^+$ and
moreover that $l_{\delta^*} \ge \rho/2>0$ as $\t \to 0$.
Thus we can apply  Proposition~\ref{prop:propA} again  to $\delta$ and
$\bch^-$ to show that
 $d(Q, \tilde \delta^*) \le c_0(\rho/2)i(\delta,\nu)\theta$ for all $Q
\in \tilde 
\delta^-$. Combining these results and observing that perpendicular
projection from $ \tilde \delta^{\pm}$ to $ \tilde \delta^*$ is
surjective, we deduce that $d(P, \tilde \delta^-) \le
c_0(\rho/2)(i(\g,\mu) + i(\g,\nu)) \t$, for all $P \in \tilde 
\delta^+$.

On the other hand, since $\mu,\nu$
fill up $\Sf$, by Lemma~\ref{lemma:intbound0} we have $i(\omega,
\mu)+i(\omega, \nu)>c_3$, so that either $i(\omega, \mu)>c_3/2$ or
$i(\omega, \nu)>c_3/2$. Assume the first
case holds.
Let    $\s$ be a segment along the lift of the axis of
$\tilde \omega^+$ of length $\rho$. Then $\s$ contains (roughly) 
$\rho/l_{\omega^+}$ periods of $\omega$, so that 
$i(\s,\mu) >
\rho i(\omega, \mu)/l_{\omega^+}$. 
  By Proposition~\ref{prop:lowerbound1}, since we also arranged $\rho <
\epsilon_2$, for any point $P \in \s$ we
have $d(P, \bch^-) \ge
 c_2i(\s, \mu) \t  $. 
 Comparison with the previous estimate applied to  
 the intersection point $P_0$
of $\omega$ and
$\delta$ on $\bch^+$ gives   
 $ l_{\omega^+}  \ge        c_2 c_3 \rho/(2(c_0(\rho/2)
(||\mu||+||\nu||))$.

Finally, if  $i(\omega, \mu)< c_2/2$ then  $i(\omega,
\nu)>c_2/2$. If  $l_{\omega^-}(p^+(\t) > \rho$ there is nothing to
prove; otherwise arguing as above
but with    $\s$ a segment along the lift of 
$\tilde \omega^-$ gives the result.
\end{proof}

The following corollary will be useful. (This can easily  be
strengthened to an assertion about the maximum distance between $\bch^+$
and $\bch^-$, but we shall not need this.)

\begin{cor} 
\label{cor:bounddistance} Suppose $P \in \bch^+$ lies on one of the
curves $\g_i^+$, $\g_i  \in \Gamma$. Then there exists a 
constant $c_4$ such that $d(P,
\bch^-) < c_4(||\mu||_{\Gamma} + ||\nu||_{\Gamma}) \t$ as $\t
\to 0$.
\end{cor}
\begin{proof} Suppose  $P \in \bch^+$ and write $\g = \g_i$. Choose
$\rho_* $ as in the above proposition. Since
$l_{\g} > \rho_* $, by
Proposition~\ref{prop:propA}  for sufficiently small $\t$
we have 
$d (P, Q) <  i(\mu, \g )O(\t)  $ for some $Q \in \tilde \g^*$, and
$l_{\g^*} >   
\rho_* /2$. Then $l_{\g^-} >\rho_* /2$. Noting that perpendicular
projection from 
 $\tilde \g^-$ to $\tilde \g^*$ is surjective,  we can apply 
Proposition~\ref{prop:propA} again to $\bch^-$ to deduce that
there is a point $P' \in \tilde \g^-$ so that  
$d (P', Q) <i(\nu, \g )O(\t) $.
Combining these results gives  $d (P, \bch^-) < (||\mu||_{\Gamma} +
||\nu||_{\Gamma})O(\t) $ as claimed.
\end{proof}

We now need to get  upper
 bounds the lengths $l_{\mu^{\pm}}$ and $l_{\nu^{\pm}}$. To do this  it
is
convenient to organise things so that the laminations $\mu$ and $\nu$
are confined to
narrow paths on the surface. As sketched by Thurston~\cite{ThuN} p.8.52,
given a geodesic lamination $\mu \in \ML$, one can find $\epsilon>0$
and a train track $\tau$, such that the $\epsilon$-neighbourhood
$N_{\epsilon}(\tau)$
is just the product of an open interval with $\tau$ (i.e. a `strip with
switches') and such that
$|\mu| \subset N_{\epsilon}(\tau)$. It is important for us to understand
the dependence of $\epsilon$ on the geometry of $\Sf$ and $\mu$, so we
state this precisely as:
\begin{prop}
\label{prop:thintrack} Suppose that $\mu \in \ML(\Sf)$ and that
$\Sf_0$ is a hyperbolic structure on the surface $\Sf$ such that
the non-cupsidal injectivity radius $inj(\Sf^C_0)> \rho$. Then there
exists
$\epsilon>0$, depending only on $\rho$
and the topology of $\Sf$, and a train track $\tau$, such that
$N_{\epsilon}(\tau)$ is  homeomorphic to $\tau$ and such that every leaf
of
$\mu$  is contained in $N_{\epsilon}(\tau)$. There is a  fixed upper
bound to the  number of switches and branches of $\tau$. 
\end{prop}
\begin{proof} The proof is explained in detail in~\cite{Penner} Theorem
1.6.5, however the dependence on ${\rm inj}(\Sf^C_0)$ is not spelled
out. 
The idea is that the complement of $|\mu|$ in $\Sf$ consists of a finite
number of ideal polygons (possibly with punctures). The area of the
$\epsilon$ neighbourhood $U_{\epsilon}$ of the boundary of any such
polygon tends to zero with $\epsilon$.

On the other hand, the bound on injectivity radius means that 
any disk $D_{\rho}$ of radius ${\rho}$ and contained in $\Sf^C_0$ is
embedded. Since
such a disk has definite area, it cannot be contained in $U_{\epsilon}$
as $\epsilon \to 0$. Thus $D_{\rho}$ must intersect $U_{\epsilon}$ in
thin
tubular  neighbourhoods of possibly branched $1$-manifolds.\footnote{
Figure
1.6.3 in~\cite{Penner} is somewhat deceptive since examination of the
constants
shows that for this to work one must take $\epsilon = O({\rho}^2)$ so
that
typically $D_{\rho}$ will contain $O(1/\sqrt \epsilon)$ such strips, not
one
as in the picture.} Now use the fact that we may choose
the neighbourhoods of the cusps such that every simple geodesic, in
particular every leaf of $|\mu|$, is entirely contained in $\Sf^C_0$. 
Clearly  the bound on $\epsilon$ depends only on ${\rho}$ and not on
$\mu$.  
The bounds on the number of branches and switches come from the obvious
bounds on the number of sides and cusps of the complementary regions
to $\mu$, which are again independent of $\mu$.
\end{proof}

We call a train track $\tau$ chosen as in the above proposition an  {\em 
$\epsilon$-thin train track}, and we say that $\mu$ is {\em carried} by
$\tau$.  It is also part of the above construction that
$N_{\epsilon}(\tau)$ is foliated by arcs of length at most
$2\epsilon$ transverse (and approximately orthogonal) to the branches.
If $b$ is a branch of $\tau$, we write $N_{\epsilon}(b)$ for the union
of the leaves transverse to  $b$.  Each of these arcs is a transversal
to $\mu$ and carries the same transverse weight which we denote
$\mu(b)$.  We emphasize that the topology of the   $\tau$
and the $\mu$ weights of each of its branches may well change as $\Sf$
varies. 

Now as in Section~\ref{sec:main}, we are seeking an upper bound on the
length of $\mu$. For fixed $\epsilon$, the branch $N_{\epsilon}(b)$ has
a definite width and thus its length is uniformly bounded above. Thus
the only way in which the lamination $\mu$ can get very long is to have
large weight concentrated in thin strips, in other words for the weights
$\mu(b)$ to get large. We rule out this possibility by  once again
playing off the
upper and lower bounds on the distance between points on $ \bch^+$ and
$\bch^-$:

\begin{prop}
\label{prop:weightbound} Choose $\rho_*$ as in
Proposition~\ref{prop:injradius}  and fix $\epsilon_*$ depending on
$\rho_*$ as in Proposition~\ref{prop:thintrack}, and so that
$2\epsilon_* <  \epsilon_2 $ as in Proposition~\ref{prop:lowerbound1}.
Suppose that $\mu \in \ML(\Sf)$ is carried on
some $\epsilon_*$-thin train track $\tau$ on the surface $\
\bch^+/G(\t)$.
Then the weight $\mu(b)$ of a branch of $\tau$ is uniformly bounded
above as $\t \to 0$.
\end{prop}
\begin{proof} As usual, let $\Gamma$ be our fixed set of curves which
fill up $\Sf$. First  assume that $N_{\epsilon_*}(b) \cap \Gamma \neq
\emptyset$.  For each $\t$, pick $P= P(\t) \in  
N_{\epsilon_*}(b) \cap \tilde \gamma_i^+$ for
some $i$.
By Proposition~\ref{prop:injradius}, the length of each $\gamma_i$ is
uniformly bounded below and so we can apply
Corollary~\ref{cor:bounddistance} to
show that $d(P, \bch^-) \le O(\t)$ with uniform constant
independent of $\t$. 
On the other hand, since the transverse measure of an arc of length
$2\epsilon_*$ containing $P$ is $\mu(b)$, by
Proposition~\ref{prop:lowerbound1} we have 
$d(P , \bch^-) > \mu(b) O(\t)$, also with a uniform constant.
Comparing these two inequalities gives a uniform upper bound on
$\mu(b)$.

Now assume that $N_{\epsilon_*}(b) \cap \Gamma = \emptyset$.
In this case $N_{\epsilon_*}(b) $ is contained in a  
 component $R$ of  $\Sf - \Gamma$ which  topologically is either a disk
or punctured disk. 
The boundary $\partial R$ consists of a bounded number of finite 
arcs $\a_j$, each contained in $\gamma_i$ for some $i$. 
We need to guard against the possibility that leaves of $\mu$ wrap
around many times inside $R$, assigning unduly large weight to
$N_{\epsilon_*}(b) $.

 Suppose first that $R$ is simply connected. Let $\lambda$ be a
connected component of the intersection of some leaf of   $|\mu|$ with
$R$. Since $\epsilon_*  $ is certainly
less than the non-cuspidal
injectivity radius, it is easy to see that $\lambda$  cannot return
within distance $\epsilon_*$ of itself.  
 From this we deduce that $\lambda$ has finite length and that it
intersects any transversal $T$ to $N_{\epsilon_*}(b) $ at most once.
It follows that there is a well defined  map  from $T \cap
|\mu|$ to $\partial R$, which sends the point $x \in T$ to the point at
which the component of $|\mu| \cap R$ through $x$ meets $\partial R$.
The map is injective when lifted to the closure of $\tilde R$
in $\HH^2$. 
Pushing forward the measure of $T$ to $\partial \tilde R$ and using the
injectivity, 
 we find $ \mu(b) \le   \sum_j i(\mu,\a_j)$. Since $\a_i \subset \g_i$
and since a  segment of $\gamma_i$ can appear at most twice on $\partial
\tilde R$, we
have  $ \mu(b) \le 2 \sum_i i(\mu,\g_i) =
2||\mu||$.

Now suppose that  $R$ contains a puncture $P$. We claim that there is no
circuit in $\tau$ encircling $P$; for otherwise by following the
boundary leaf of $|\mu|$ nearest to $P$ round this circuit we would have
a leaf of $|\mu|$ parallel to a loop around $P$, which is impossible.
Thus we can find an arc $\a$ joining $P$ to $\partial R$ with   $\Int \a
\subset R$ and   $\a \cap |\mu| = \emptyset$. Now  repeat the
above argument working in $R - \a$. 
\end{proof}

\begin{cor}
\label{cor:lengthbound1} Suppose that $\mu,\nu \in \ML$ and that $G(\t)=
G(\t \mu, \t
\nu) \in \QF(\Sf)$. Then  there is a uniform upper bound to the lengths
$l_{\mu^+} = l_{\mu^+}(p^+(\t))$ as $\t \to 0$.
\end{cor}
\begin{proof} Fix $\epsilon_*$ as in Proposition~\ref{prop:thintrack} 
and find an $\epsilon_*$-thin train track $\tau$ on the surface
$p^+(\t)$ carrying $\mu$.
The length of $\mu$ on $p^+(\t)$  is clearly estimated from
above by 
$\sum_i \mu(b_i) l_{b_i}$ where $l_{b_i}= l_{b_i}(p^+(\t))$ is the
length of the branch
$b_i$ and $\mu(b_i)= \mu(b_i)(p^+(\t)) $ its weight.  Now the strip of
width $2\epsilon_*$ about $b_i$ is embedded on
$\Sf$, so since $\epsilon_*$ is fixed independent of $\t$ the length
$l_{b_i}$  is
uniformly bounded above by an area estimate. The result follows.
\end{proof}

\noindent {\sc Proof of  Theorem~\ref{thm:main}.}
By Corollary~\ref{cor:lengthbound1}, the lengths 
 $l_{\mu^+}(\t)  $ and $l_{\nu^-} (\t)  $ are uniformly bounded above.
 We claim that this implies that the hyperbolic structures $p^{\pm}(\t)
$ of $\bch^{\pm}/G(\t)$ both converge in $\F$, and that their limits
coincide.

Suppose that the sequence $p^+(\t) $ does not converge, then (after
passing to a subsequence if necessary)  it limits on  some projective
lamination $[\xi] \in \PML$. This means that there exists a lamination
$\xi \in [\xi]$ and a sequence $h_n \to 0$ such that $h_n l_{\g^+}(\t)
\to i(\g,\xi)$ for all $ \g \in  \S$. 
Now it follows from  Propositions~\ref{prop:propA}
and~\ref{prop:injradius} that  $l_{\g^-}(\t) /l_{\g^+}(\t)  \to 1$ as
$\t \to 0$. Thus  $h_n l_{\g^-}(\t)  =  h_n l_{\g^+}  \cdot l_{\g^-}
/l_{\g^+}  \to
i(\g,\xi)$ and so $p^-(\t) $ converges to  $[\xi] \in \PML$ also. 

From  the definition of convergence to $\PML$, the fact that 
$l_{\mu^+}(\t)  $ remains  bounded   implies that  $i(\mu, \xi) = 0$;
 likewise  $i(\nu, \xi) = 0$. However since $\mu,\nu$ fill up $\Sf$,
this is impossible. 
We conclude that $p^+(\t) $  converges to a point $  p^+_{\infty} \in
\F$,
and likewise $p^-(\t) \to p^-_{\infty} \in \F$. Finally the same length
comparison
$l_{\g^-}/l_{\g^+} \to 1$ shows that  $  p^+_{\infty} = p^-_{\infty}$,
and our claim follows.

We can now use the lower bound on lengths from
Proposition~\ref{prop:injradius} to conclude the proof exactly as 
we did at the end of Section~\ref{sec:main}.
\qed

 \section{Diagonal limits}
\label{sec:diagonal}

In this final section we discuss the question of diagonal limits.
Here is a simple example which shows that care is needed.
Let $\{\a_1,\ldots, \a_k \}$ and $\{\b_1,\ldots, \b_l \}$ be two systems
of disjoint curves which fill up $\Sf$, but so that $\{\a_2,\ldots, \a_k
,\b_1,\ldots, \b_l \}$  do not. 
Fix coefficients $a_2,\ldots, a_k ,b_1,\ldots, b_l > 0$ and choose a
sequence $h_n \to 0$. Define 
$\mu_n  = h_n \a_1 + \sum_{i=2}^k a_i \a_i$ and $\nu_n \equiv \nu = 
\sum_{i=1}^l b_i \b_i$.
Obviously, $\mu_n  \to  \mu_{\infty} =  \sum_{i=2}^k a_i \a_i$ in
$\ML(\Sf)$.

Now consider the sequence of groups $G_n = G(\t \mu_n, \t \nu)$ as $\t
\to 0$. 
Since $\mu_{\infty}, \nu$ do not fill up $\Sf$,  the length function
$l_{\mu_{\infty}} + l_{\nu}$ does not have a minimum on $\Teich(\Sf)$
(see~\cite{KerckLM} p.194),   and so the sequence
$G_n$ cannot have a Fuchsian limit. In fact, because the ratio
between the weights on $\a_1$ and each of the curves in $|\nu|$ tends to
zero, the constants in the estimate for the lower bound on distance
between  a bending line  and the opposite side of $\bch$ in
Proposition~\ref{prop:lowerbound} become arbitrarily small. Thus the
length bound in Proposition~\ref{prop:lengthupperbound} fails, in other
words, $l_{\a_1} \to \infty$ 
and  the groups $G_n$ diverge. (For the limiting behaviour along lines
of minima, see~\cite{DSends}.)

The final important  step in our proof of Theorem~\ref{thm:main} was
establishing  the
length bounds on $l_{\mu^+}$ and $l_{\nu^-}$, from which we
deduced that the corresponding surfaces lay in a compact set in
$\Teich(\Sf)$.  However our example 
shows that  to prove Theorem~\ref{thm:diaglimit}, it is no use just
establishing uniform upper bounds  on the lengths $l_{\mu_n}(\t_n)$ and
$l_{\nu_n}(\t_n)$. In fact it is easy to adjust the coefficients $h_n$ 
so as to produce a sequence  $p_n \in \F$
such that  the lengths $l_{\mu_n}(p_n)$ and $l_{\nu}(p_n)$ are
uniformly bounded above,  but such that $l_{\a_1}(p_n) \to \infty$, so
that the sequence $p_n$  exits every compact set in $\Teich(\Sf)$.
To resolve this we need to do more work to get a uniform upper bound on
the lengths of the curves $\Gamma$.
We shall invoke the convergence of the
laminations through the following improved version of
Lemma~\ref{lemma:intbound0}, which allows us to avoid the
technical difficulty  that it is not clear how to control
the behaviour of transversals to $\mu$ as we transfer from surface to
surface --  {\em a priori} a narrow strip containing heavy weight on one
structure might become extremely wide on another.
We keep the same name for  the  constant $c_3$ although the actual value
may have changed.

\begin{lemma}
\label{lemma:intbound} Let   $\mu$ and $\nu$ fill up $\Sf$ and suppose
that $\mu_n \to \mu, \nu_n \to \nu$ in $\ML$. Then there exists
$c_3>0$ such that $i(\g ,\mu_n) +i(\g ,\nu_n) > c_3$ for all $\g \in \S$
and all $n$.
\end{lemma}
\begin{proof} The proof is almost the same as the previous version. Once
again we can work  entirely with a fixed hyperbolic structure $p_0 \in
\F$ whose
non-cuspidal injectivity radius is the fixed value $\rho_0>0$. 
If the result is false, then we can find a  sequence $\gamma_n \in \S$
such
that $i(\g_n,\mu_n) +i(\g_n,\nu_n) \to 0$. As before, passing to a
further subsequence we find $h_n >0 $ such that $h_n
\delta_{\g_n}  \to  \xi $ in $\ML$ and $h_n
\le 2/\rho_0 $. Then $i(h_n\g_n,\mu_n) +i(h_n\g_n,\nu_n) \to 0$
but also $i(h_n\g_n,\mu_n) +i(h_n\g_n,\nu_n) \to i(\xi,\mu)
+i(\xi,\nu)$. Since $\mu$ and $\nu$ fill up $\Sf$, 
this is impossible. 
\end{proof}

\medskip

\noindent {\sc Proof of Theorem~\ref{thm:diaglimit}.}
We have to study the behaviour of the groups $G(\theta_n \mu_n, \theta_n
\nu_n)$ as $\mu_n \to \mu$, $\nu_n \to \nu$, and $\t_n \to 0$.  Observe
that the constants involved in the proof of
Theorem~\ref{thm:main} depended only on
the topology of $\Sf$ 
and the intersection numbers $i(\mu,\g)$ and $i(\nu,\g)$
for $\g $ in the fixed finite set $\Gamma$. 
In particular, inspection of the proof of
Proposition~\ref{prop:injradius} shows that the lower bound $\rho_*$ on
the non-cuspidal injectivity radius 
depends on the universal Margulis constant, upper bounds for the norm
$||\mu||_{\Gamma}+ ||\nu||_{\Gamma}$, and the constant $c_3$ of
Lemma~\ref{lemma:intbound0}. 
If $\mu_n, \nu_n \to \mu,\nu$ in $\ML$, then 
$i(\mu_n,\d) \to i(\mu ,\d)  $ and $i(\nu_n,\d) \to i(\nu ,\d)  $
for all $\d \in \pi_1(\Sf)$.
Since $\Gamma$ is a fixed finite set, and since we have just shown that
the constant $c_3$ is independent of $n$, we deduce that all bounds in
question are uniform as $n \to \infty$.
In particular the  lower bound $\rho_*$ can be chosen  uniform for
all the  structures $p^{\pm}(\t_n)$ on the surfaces
$\bch^{\pm}/G(\theta_n \mu_n, \theta_n \nu_n)$. 

 From the results of~\cite{KerckLM}, we know that 
$M(\mu_n,\nu_n) \to M(\mu,\nu)$. 
Thus to prove diagonal convergence, we just need a uniform estimate on
the
convergence of $G(\theta_n \mu_n, \theta_n \nu_n)$  to $M(\mu_n,\nu_n)$.
From  Proposition~\ref{prop:propA} and the uniform lower bound $\rho_*$,
we get
$$| l_{\delta} (p^+(\t_n)) / l_{\delta} (p^-(\t_n))| < 1+ O(\t_n^2) $$
for any curve $\delta \in \pi_1(\Sf)$, with uniform constants depending
only on $\delta$,$\mu$ and $\nu$.

We shall show in Proposition~\ref{prop:diagbound} below that the lengths 
of any fixed curve $\delta \in \S$ on either surface $p^{\pm}(\t_n)$
have a uniform upper bound as $n \to \infty$. 
Thus we obtain $$| l_{\delta} (p^+(\t_n)) - l_{\delta} (p^-(\t_n))| <  
O(\t_n^2), $$  with constants depending on $\delta$. 
Taking a large enough finite set of  curves  to determine  the analytic
structure  on $\QF$ establishes uniformity of convergence, and 
the result follows.
 \qed

In view of the above, the main work remaining work in proving 
Theorem~\ref{thm:diaglimit} is establishing the following:
\begin{prop}
\label{prop:diagbound}  Suppose that $\mu,\nu$ fill up $\Sf$. Let $\g
\in \S$ be fixed and let $\mu_n, \nu_n \to \mu,\nu$ in $\ML$ be such
that $\mu_n,\nu_n$ fill up $\Sf$ for all $n$. Suppose that  $\t_n \to
0$. Then the lengths
$l_{\g} (p^{\pm} (\t_n)) $ are uniformly bounded above as $n \to \infty$
(with a bound depending only on $\g$, $\mu$ and $\nu$ and the
topological type of $\Sf$).
\end{prop}

As indicated by our counter example, convergence may fail if the curve
$\g$ only meets $\mu_n$ or $\nu_n$ branches of vanishingly small weight.
In fact a closed loop in  $\mu_n$ of very small weight may itself become
extremely long; this is only avoided by the hypothesis that the limit
laminations $\mu$ and $\nu$ themselves fill up the surface. Another
possibility is that $\g$ might meet a loop of definite weight and
bounded length, but by wrapping around it many times it could itself
become extremely long. Thus to prove Proposition~\ref {prop:diagbound} 
we show  first (Proposition~\ref{prop:fixeddist}), that  every segment
of $\g$ of some definite length must meet some bending line of definite
weight, and second (Proposition~\ref{prop:diagbound0}), that it must 
meet the relevant bending line sufficiently transversally  to make a
definite contribution to the positive number 
$i(\mu_n , \g)+ i(\nu_n,\g)$.

\begin{prop}
\label{prop:fixeddist}  Let     $\mu_n, \nu_n, \mu,\nu$ be as above. Fix
$\e_*$ as in Proposition~\ref{prop:thintrack} and for each $n$, suppose    
$\tau(\mu_n),\tau(\nu_n)$  are $\e_*$-thin train tracks on $p^+(\t_n)$
carrying  $\mu_n$ and  $\nu_n$ respectively.
Then there exist uniform constants  $L_1, k_1>0$ such that if $\s$ is
any geodesic segment 
on $p^+(\t_n)$  contained in a complete simple geodesic and of length
at least $L_1$, then there is  a point $P \in \s$ such that $P \in
|\mu^+_n| \cup |\nu^+_n|$ and such that $P$ is contained in a
$\tau(\mu_n)$
or $\tau(\nu_n)$ branch of transverse weight at least $k_1$.
\end{prop} 
\begin{proof}  
In the statement  $\mu_n,\nu_n$ of course refer to the representatives
of these laminations on $p^+(\t_n)$. As already observed, the number of
branches of $\tau(\mu_n)$ and
$\tau(\nu_n)$ has a uniform upper bound depending only on the topology
of $\Sf$. Moreover it is easy to see that any component of
$N_{\e}(\tau(\mu_n)) \cap N_{\e}(\tau(\nu_n))$ contains a ball of radius
at least $\e$, so that  by an area argument the total number of
intersection points of $\tau(\mu_n)$ and $\tau(\nu_n)$  has an upper 
bound independent of $n$.
It follows that  there is also a uniform upper bound to the number of
sides
of each complementary region of $\tau(\mu_n) \cup \tau(\nu_n)$.
As in the proof of Corollary~\ref{cor:lengthbound1}, by area
considerations there is a uniform upper bound $l_0$ to the length of any
branch of $\tau(\mu_n)$ or $\tau(\nu_n)$. Now each complementary region
is either simply connected or a once punctured disk. Thus every simple
arc crossing a complementary region is homotopic to a path along the
boundary but not fully encircling the boundary. (No simple
geodesic can completely encircle the puncture.) 
This gives a uniform upper bound $l_1$ to the length of the intersection
of a simple geodesic with each complementary region; we may as well
assume that $l_1 > l_0$.

Now let $\s:[0,T] \to \bch^+(\t_n)$ be a geodesic segment parameterized
for
convenience by arc length.  We associate a crude symbol sequence to $\s$
as follows. First, after removing a transverse arc through each switch
of
$\mu$, the open neighbourhood $N_{\e}(\tau(\mu_n))$ is disconnected into
a finite number of open sets $V_i$,  one for each branch of $\mu$.
Disconnect $N_{\e}(\tau(\nu_n))$ into sets $W_j$  a similar way and let
  $ \cal B$ denote the set of all components of the resulting dissection
of $N_{\e}(\tau(\mu_n)) \cup N_{\e}(\tau(\nu_n))$; thus a set in $\B$ is
a component either of $V_i \setminus N_{\e}(\tau(\nu_n))$, or of $W_j
\setminus N_{\e}(\tau(\mu_n))$; or of $V_i \cap  W_j$.   As we have seen
the size of  $ \cal B$ is uniformly bounded above by some $ M \in \NN$.  
The  points at which  $\s$ meet $\partial Y$ for any  $Y \in \cal B$
give a partition of $\s$  at the points $0= t_0 < t_1< \ldots < t_m=T$. 
Thus each open arc $\s(t_i , t_{i+1})$ is contained either in a
component  $Y \in \B$, or in a complementary region of $\tau_{\mu_n}
\cup
\tau_{\nu_n}$.
We associate to $\s$ the sequence $e_1 \ldots e_m$ where $e_i =Y$ if 
$(t_{i-1},t_{i}) \subset Y$ and $e_i =X$ if $\s(t_i,t_{i+1})$ is
contained  in a complementary region.  Observe that from the definition,
no symbol is immediately followed by itself.
Also note that $T= l_{\s}  \le   ml_1$.

Now consider any simple geodesic segment $\s$ with length $l_{\s} \ge
(2M+1)l_1$.
Suppose its symbol  sequence is  $e_1  \ldots  e_m$. Then $(2M+1)l_1 \le
l_{\s} \le  m l_1 $ so that $M+1 \ge [m/2]$. Since the symbol $X$  never
follows itself, at least $[m/2]$  symbols from $e_1  \ldots  e_m$ belong
to $\B$ and hence  
 some symbol $b  \in \B$ occurs twice; that is, some subarc $\s' \subset
\s$ runs from the component  $Y$ to itself. Assume that $\s'$ is a
minimal segment of this type, in the sense that the length of its symbol
sequence is least possible, so that in particular  this length is at
most $2M+1$.   Let $\s_{Y}$ be the geodesic arc joining
the first point at which $\s'$ leaves  $\partial Y$ to the next point at
which it reenters it. Since $Y$ is geodesically convex, $\s_Y \subset
Y$. Thus $\s'' = \s' \cup
\s_Y$ is a loop of length at most  $(2M+2) l_1$.

By Lemma~\ref{lemma:intbound}, we have $i(\s'',\mu_n)+ i(\s'',\nu_n)>
c_3 $.  It is clear that
$i(\s'',\mu_n) = i(\s',\mu_n) + i(  \s_Y,\mu_n)$
and similarly for $\nu_n$.  
Thus either $i(\s_Y,\mu_n) +i(\s_Y,\nu_n) > c_3/2$ or $i(\s',\mu_n)+
i(\s',\nu_n)> c_3/2$.
In the first case we see that $Y$ lies in a branch of either 
$\mu_n$ or $\nu_n$ of weight at least $c_3/4$.
In the second case, since $\s'$ contains at most $(2M+1)$ elements in
its symbol sequence, we see that it must meet at least one branch $b'
\in \B$ of $\mu_n$ or $\nu_n$-weight at least  $c_3/2(2M+1)$.
Thus in all cases $\s'$ contains some point   which lies in a $\mu_n$ or
$\nu_n$  branch  of weight at least $ c_3/2(2M+1)$. Setting
$L_1 = (2M+1)l_1$ and $ k_1= c_3/2(2M+1)$ gives the result.
\end{proof}

In the proof of the next proposition we shall need to transfer the
lamination $\nu$ from $\bch^+$ to  $\bch^-$.
The following shows that we can do this without serious loss of control.
For clarity, we denote the copies of a   laminations  $\xi$  on
$\bch^{\pm}$ by $\xi^{\pm}$ respectively.

\begin{lemma}
\label{lemma:transfer} Let $\lambda^+$ be a leaf of the lamination
$\nu^+_n$ lifted to the surface $\bch^+(\t_n) $ and  
let $\lambda^-$ be  the corresponding leaf 
 of $\nu^-_n$ on the surface $\bch^- (\t_n)$, so that $\lambda^{\pm}$
have
the same endpoints on $\partial \HH^3$. Then for any point $P \in \l^+$
we have 
 $d(P,\lambda^-) \le c_5\t_n$ with a uniform constant $c_5$ as $n \to
\infty$.
 \end{lemma}
\begin{proof}
The idea is obviously to imitate the proof of
Proposition~\ref{prop:propA}. To do this we need to see that 
there is a uniform upper
bound to the $\mu^+_n$ mass of any geodesic segment on $\bch^+(\t_n)$ of
definite length at most $1$ say.
In fact if $\tau$ is an $\epsilon_*$-thin train track  carrying
$\mu^+_n$,
by Proposition~\ref{prop:weightbound}
there is a uniform upper bound to the weight of each $\mu^+_n$ branch,
moreover
away from $\epsilon_*$ balls around the  switches  there is  
by the construction of $\tau$  a uniform lower bound to the distance
between leaves of
$\mu^+_n$ contained in distinct branches. This gives the required bound.  

Now pick equally spaced points $P_m, m\in
\ZZ$  at unit distance apart along $\l^+$. The above discussion gives a
uniform upper bound to $i(\s_m,\mu_n)$ where $\s_m$ is the segment of
$\l^+$ from $P_m$ to $P_{m+1}$. Thus we may  argue exactly as in the
proof
of Proposition~\ref{prop:propA} to show that 
all points on $\l^+$ are at most a uniform distance $O(\t)$ away from
the corresponding leaf $\l^-$.  
\end{proof}

\begin{prop}
\label{prop:diagbound0} Let $\g \in \S$. Then there exist  $L_0, C_0>0$,
depending only on $i(\g,\mu), i(\g,\nu),i(\nu,\mu)$, such that if $\s$
is any geodesic segment contained in $\g^+$ of length at least $L_0$ on
any of the  hyperbolic surfaces $p^{\pm}(\t_n) $, then $i(\s,\mu_n)+
i(\s,\nu_n)>C_0$.
\end{prop}
\begin{proof} 
If the result is false, then we can find  a structure $p^+(\t_n)$ say on
which $\g$ has an arbitrarily  long segment $\s$ for which $i(\s,\mu_n)
+ i(\s,\nu_n)$ is arbitrarily small. The argument will follow  the same
lines as that of Proposition~\ref{prop:lengthupperbound}.

Consider such a segment $\s$ where the choice of constants will be
determined later,  and let $\hat \s$ be the $\HH^3$ geodesic joining its 
endpoints $X$ and $X'$. Clearly we may assume that $L_0>1$ say; then as
in the first part of the proof of Proposition~\ref{prop:propA}, there is
a universal constant $c$ such that $d(P, \hat \s) < ci(\s,\mu_n) \t$ for
all $P \in \s$, and such that $l_{\hat \s} > (1-c(i(\s,\mu_n)^2\t^2)
l_{\s} > l_{\s}/2$ say for all small enough $\t$. 
Choose $L_1$ as in Proposition~\ref{prop:fixeddist}. 
 By Lemma~\ref{lemma:cegdistance}, given $h>0$  we can find $L_2 =
L_2(h)$
such that, if $l_{\s} >L_2$, then  $d(Q,\tilde \g^*) < h \t$ for all
points $Q$
on a segment $\hat
\s' \subset \hat \s$ of length at least  $L_1$. Since  perpendicular
projection from $\s$ to $\hat \s$ is surjective,  we can find  a
subarc of $\s_1$ of $\s$ length at least $L_1$ for which 
$d(Q, \tilde \g^*) < (h + ci(\s,\mu_n)) \t$  for all $Q \in \s_1$.
Let $Y,Y'$ be points in $\bch^-$ close to $X,X'$; arguing similarly we
can find a point $R$ in the
$\bch^-$ arc from $Y$ to $Y'$ such that   $d(Q, R) < (2h +
c(i(\s,\mu_n)+i(\s,\nu_n))) \t$.

By Proposition~\ref{prop:fixeddist}, the arc $\s_1$ contains a point $P$
which lies within distance at most $\e_*$ of a point in a branch of
either
$\mu^+_n$ or $\nu^+_n$ and of weight at least $k_1$.
Suppose first this is a $\mu^+_n$-branch. Then by
Proposition~\ref{prop:lowerbound1}, there is a constant $c_2>0$  such
that $d(P,\bch^-) > c_2 k_1 \t $.
Choose $h$ as above with $2h < c_2 k_1/2$ and  then use our hypothesis
to choose
$\s$   with $l_{\s} > L_2(h)$ and  $i(\s,\mu)+i(\s,\nu)< c_2 k_1/2$. 
Comparison with the estimate $d(Q, R) < (2h +
c(i(\s,\mu_n)+i(\s,\nu_n))) \t$ gives a contradiction.

 If $P$ is in a branch of $\nu^+_n$ of weight at least $k_1$, we can
clearly make a similar argument provided we can find a point $P' \in
\bch^- $ near $P$ and within distance $  \e_2$ of a $\nu^-_n $ branch of
definite weight.  From Lemma~\ref{lemma:transfer}  we
see that  we can find
a transversal $T'$ to $ \nu^-_n$ of $\bch^-$ length at
most $2\e_* +2c_5 \t $ and with $\nu^-_n(T')>k_1$. Provided $\t $ is
sufficiently small we have the result. 
\end{proof}

\noindent {\sc Proof of Proposition~\ref{prop:diagbound}.} 
From Proposition~\ref{prop:diagbound0}
we easily obtain the bound
$l_{\g^+}(p^+(\t_n) ) \le 2 (i(\g,\mu)+ i(\g,\nu))L_0/C_0$.  
\qed
 \noindent This completes the proof of
Theorem~\ref{thm:diaglimit}.

\small{
}\end{document}